\documentclass[12pt,a4paper]{article}
\usepackage{amssymb,amsmath}
\usepackage{amsfonts}
\usepackage{amsthm}
\usepackage{a4}
\usepackage{tikz}
\usepackage{enumerate}%
\usepackage{verbatim}
\usepackage{mathalfa}
\usepackage{caption}
\usepackage{float}
\newcommand{\para}{\par\vspace{.25cm}}

\newtheorem*{theorem*}{Theorem}
\newtheorem{theorem}{Theorem}
\newtheorem{lemma}{Lemma}
\newtheorem{cor}{Corollary}
\newtheorem{remark}{Remark}

\newcommand{\Q}{\mathbb{Q}}

\newcommand{\C}{\operatorname{Cen}}

\newcommand{\G}{\operatorname{Gal}}

\begin{document}
	\title{\bf A computational approach to Brauer Witt theorem using Shoda pair theory}
	\author{ Gurmeet K. Bakshi{\footnote {Research supported by Science and Engineering Research Board (SERB), DST, Govt. of India under the scheme Mathematical Research Impact Centric Support (sanction order no MTR/2019/001342) is gratefully acknowledged.}}  and Jyoti{\footnote {Research supported by Council of Scientific and Industrial Research (CSIR), Govt. Of India under the reference no. 09/135(0886)/2019-EMR-I is gratefully acknowledged.} \footnote{Corresponding author}} \\ {\em \small Department of 
			Mathematics,}\\
		{\em \small Panjab University, Chandigarh 160014, India}\\{\em
			\small email: gkbakshi@pu.ac.in and jyotigarg0811@gmail.com} }
	\date{}
	{\maketitle}
	
\begin{abstract} \noindent
A classical theorem due to Brauer and Witt implies that every simple component of the rational group algebra $\mathbb{Q}G$ of a finite group $G$ is Brauer equivalent to a cyclotomic algebra containing $\mathbb{Q}$ in its centre. The precise description of this cyclotomic algebra is not available from the proof of the Brauer-Witt theorem and it has been a problem of interest to determine the same in view of its central role in the study of central simple algebras. In this paper, an approach using Shoda pair theory is described, which is quite efficient from computational perspective. 
\end{abstract}	
\para \noindent {\bf MSC2000 :} 16S34, 16K20, 16S35, 20C05, 16K50.
\section{Introduction}
Let $G$ be a finite group and $\Q G$ the rational group algebra of $G$. Given an irreducible character $\chi$ of $G$, it is well known (\cite{TY}, Proposition 1.1) that  $$e_{\mathbb{Q}}(\chi) = \frac{\chi(1)}{|G|}\sum_{\sigma \in \operatorname{Gal}(\mathbb{Q}(\chi)/ \mathbb{Q})} {\sum _{g \in G}} \sigma(\chi(g))     g^{-1}$$ is a  primitive central idempotent of $\Q G$, called the primitive central idempotent of $\Q G$  realized by $\chi$. Here $|G|$ is the order of $G$, $\Q(\chi)$ is the field obtained by adjoining the character values $\chi(g),~g \in G,$ to $\Q$ and $\G(\Q(\chi)/\Q)$ is the Galois group of extension $\Q(\chi)$ over $\Q.$
A classical result due to Brauer and Witt (\cite{Jes}, Theorem 3.7.1) implies that the simple component $\mathbb{Q}Ge_{\mathbb{Q}}(\chi)$ of $\Q G$ is Brauer equivalent to a cyclotomic algebra containing $\Q$ in its center.  
The proof of Brauer-Witt theorem does not provide any method to determine the underlying cyclotomic algebra or even its order in the Brauer group, i.e., its Schur index (denoted by $m_{\Q}(\chi)$) in terms of the group structure of $G$. One of the crucial step in the  proof of Brauer-Witt theorem is  the existence of a $p$-elementary section of the group $G$ that determines the $p$-part of $m_{\mathbb{Q}}(\chi)$  for each prime $p$ dividing it. However there is no method to determine this $p$-elementary subgroup in order to understand the cyclotomic algebra which is involved. A slightly different approach was later given by Olteanu \cite{Gab2007} in which instead of using $p$-elementary subgroup, strongly monomial characters of subgroups of $G$ are used. The reason for this shift from $p$-elementary subgroups to strongly monomial character of subgroups is due to a good understanding of the simple component of the rational group algebra corresponding to a strongly monomial character from the work of Olivieri, del R{\'{\i}}o and Sim{\'o}n \cite{OdRS04}. The approach of Olteanu heavily depends on the search with GAP in order to find the strongly monomial characters of subgroups. \para \noindent The main result of this paper is a generalization of the work of Olivieri, del R{\'{\i}}o and Sim{\'o}n (\cite{OdRS04}, Proposition 3.4) and thus provide a computational approach to Brauer-Witt theorem for a substantially large class of monomial groups. We have proved in section 3 that the simple component of the rational group algebra of a generalized strongly monomial character is isomorphic to a matrix algebra over a cyclotomic algebra and a  description of the cyclotomic algebra from the subgroup structure of $G$ has also been provided very precisely. It needs to be pointed out that the generalized strongly monomial characters were defined by the first author with Kaur in \cite{BK3} and they proved in \cite{BK4} that the generalized strongly monomial groups, i.e., monomial groups where all irreducible characters are generalized strongly monomial, is quite a vast  class of monomial groups. For example, in \cite{BK4}, it is proved that all well known classes of monomial groups like nilpotent-by-supersolvable, strongly monomial, subgroup closed monomial, subnormally monomial, supermonomial, monomial groups with Sylow towers and many more are  generalized strongly monomial. Also from the properties of generalized strongly monomial groups given in \cite{BK4}, it seems quite hard to find a monomial group  which is not generalized strongly monomial, thus justifying that the given generalization is significant.\para \noindent In general, it is hard to compute the Schur index $m_{\Q}(\chi)$ of a character $\chi$ for a given group $G$ and Brauer-Witt theorem is one of the tools to handle the same. For some works in this direction, see (\cite{Herman1995MetabelianGA}, \cite{herman_1996}, \cite{Gab2007}, \cite{shi}, \cite{Spiegel1980OnTS}, \cite{Unger2019AnAF}, \cite{Yamada1978MoreOT}, \cite{yamada1978schur}). In section 4, we consider a semidirect product of the extraspecial $p$-group of order $p^3$ with the cyclic group of order $2^n$, where $p$ is an odd prime congruent to $1$ modulo $4$ and ${n-1}$ is the highest power of $2$ dividing $p-1$. We use the theory developed in section 3 to precisely write each simple component of its rational group algebra as a matrix algebra over a cyclotomic algebra. This description allows us to tell the Schur indices of all of these simple components and provides the complete Wedderburn decomposition of the rational group algebra of this  group of order $p^32^n$. It may be pointed out that the smallest known example of a monomial group which is not strongly monomial, i.e., SmallGroup$(1000,86)$, is a  particular case of this example when $p=5$. \para \noindent In section 5, we consider a semidirect product of $Q_8 \times C_7$ with $C_3$, where $Q_8$ is the quaternion group of order $8$ and $C_n$ is the cyclic group of order $n$ for any $n \geq 1$. In \cite{Gab2007}, Olteanu computed the simple component of the rational group algebra of this group realized by a monomial character which is not strongly monomial by finding strongly monomial characters of subgroups using a GAP search. We have revisited this example here and have done the precise computation theoretically without any GAP search, thus illustrating the theory developed in this paper. \para \noindent
\section{Preliminaries}
Throughout this paper, $G$  denotes a finite group.  A \textit{monomial character} $\chi$ of $G$ is the one which is induced from a linear character of a subgroup of $G$. In case $\chi = \lambda^G$, where $\lambda$ is a linear character on a subgroup $H$ of $G$, then $H$ and $K= \operatorname{ker} \lambda,$ the kernel of $\lambda,$ satisfy the following:
\begin{itemize}
	\item [(i)] $K \unlhd H$,  $H/K$ is cyclic;
	\item [(ii)] if $g \in G$ and $[H, g] \cap H \subseteq K,$ then $ g \in H.$ Here $[H,g]=\langle g^{-1}h^{-1}gh, h \in H \rangle.$\end{itemize}
A  pair $(H,K)$ of subgroups of $G$ which satisfies (i) and (ii) above is called a \textit{Shoda pair} of $G$ (\cite{OdRS04}, Definition 1.4) and we call it a Shoda pair of $G$ arising from $\chi$. This gives us an understanding that the Shoda pairs arise from the monomial irreducible characters of $G$. Conversely, if  $(H,K)$ is a Shoda pair of $G$, then by (\cite{CR}, Corollary 45.4), we have $\lambda^G$ is irreducible for any linear character $\lambda$ on $H$ with kernel $K$ and we call it an irreducible character of $G$ arising from the Shoda pair $(H,K).$ 
\para \noindent  Suppose $(H,K)$ is a Shoda pair of $G$ and $\lambda$, $\lambda'$ are two linear characters on $H$ with kernel $K$, then clearly $\lambda'= \sigma\circ \lambda$ for some automorphism $\sigma$ of $\Q(\lambda)$ and thus $e_{\mathbb{Q}}(\lambda^G)  = e_{\mathbb{Q}}(\lambda'^G).$  Hence all the irreducible characters of $G$ which arise from the Shoda pair $(H,K)$ give us the same primitive central idempotent of $\mathbb{Q}G$, which we call as the primitive central idempotent of $\mathbb{Q}G$ realized by the Shoda pair $(H,K)$. Furthermore, we call the corresponding simple component $\Q Ge_{\Q}(\lambda^G)$ as the simple component of $\Q G$ realized by the Shoda pair $(H,K)$.\para \noindent 
Two Shoda pairs  of $G$ are said to be \textit{equivalent} if they realize the same primitive central idempotent of $\mathbb{Q}G$.  The following is a necessary and sufficient condition on the equivalence of Shoda pairs (\cite{Jes}, Problem 3.4.3): \begin{quote} 
	\textit{If $(H_{1},K_{1})$, $(H_{2},K_{2})$ are Shoda pairs of $G.$ Then $(H_{1},K_{1})$ is equivalent to $(H_{2},K_{2})$ if and only if $H_{1}^g \cap K_{2}= K_{1}^g \cap H_{2}$, for some $g \in G,$ where $H^g=g^{-1}Hg$.} 
\end{quote}
\begin{remark}{\label{remark1}}
	(a) If $H$ is a normal subgroup of $G$ and $(H,K_{1})$, $(H,K_{2})$ are Shoda pairs of $G$, then they are equivalent if and only if $K_{1}$ and $K_{2}$ are conjugate in $G$.\\
	(b) If $(H_i,K_i),~i=1,2,$ are two equivalent Shoda pairs of G and $\lambda_{i}$ is a linear character on $H_{i}$ with kernel $K_{i},$ then $\lambda_{1}^G$ and $\lambda_{2}^G$ have the same character degree. This is because $e_{\Q}(\lambda_{1}^G)= e_{\Q}(\lambda_{2}^G)$ implies that $\lambda_{1}^G=\sigma \circ \lambda_{2}^G$ for some automorphism $\sigma$ of $\Q(\lambda_{2}^G).$ \end{remark}
\noindent By a {\it complete and irredundant set} of Shoda pairs of $G$, we mean  a set of representatives of distinct equivalence classes of  Shoda pairs of $G$. If $G$ is a monomial group and $\mathcal{S}$ is a complete and irredundant set of Shoda pairs of $G$, then $$ \mathbb{Q}G \cong \oplus_{(H,K)\in \mathcal{S}} \mathbb{Q}G e_{\mathbb{Q}}(\lambda^G)$$ is the Wedderburn decomposition of $\mathbb{Q}G $, where for any $(H,K) \in \mathcal{S}$, $\lambda$ denotes a linear character on $H$ with kernel $K$. \para \noindent For a ring $R$, let $M_{n}(R)$ denote  the ring of $n\times n$ matrices over $R$. Also, let $\operatorname{Aut}(R)$  and $\mathcal{U}(R)$ denote the group of automorphisms of $R$ and the group of units of $R$ respectively. Let $\sigma: G \rightarrow \operatorname{Aut}(R)$ and $\tau : G \times G \rightarrow \mathcal{U}(R)$ be two maps which satisfy the following relations,
$$\tau_{gh,x}\sigma_{g}(\tau_{g,h})= \tau_{g,hx} \tau_{h,x}$$ and $$\tau_{g,h}\sigma_{g}(\sigma_{h}(r))= \sigma_{gh}(r) \tau_{g,h},$$for all $g,h,x \in G$ and $r \in R$. Here $\sigma_{g}$ is the image of $g$ under the map $\sigma$ and $\tau_{g,h}$ is the image of $(g,h)$ under the map $\tau.$ Let $R*_{\tau}^{\sigma} G$ denote the set of finite formal sums $\{ \sum z_{g}a_{g} ~|~ a_{g} \in R, g \in G \},$ where $z_{g}$ is a symbol corresponding to $g \in G$. Equality and addition in $R*_{\tau}^{\sigma} G$ are defined componentwise. For $g,h \in G$ and $r \in R,$ by setting
$$z_g z_h = z_{gh} \tau_{g,h},$$ $$rz_g = z_g \sigma_g(r)$$ and extending this rule distributively, $R*_{\tau}^{\sigma}G$ becomes an associative ring, called the {\it crossed product} of $G$ over $R$ with twisting $\tau$ and action $\sigma$. We call $\{z_{g}~|~g \in G\}$ as a set of basis units of $R*_{\tau}^{\sigma} G$.\para \noindent For any ring $R$, denote by $\mathcal{Z}(R)$ the center of $R$. In general, given a crossed product $R*^{\sigma}_{\tau}G,$ $\mathcal{Z}(R*^{\sigma}_{\tau}G)$ contains $\mathcal{Z}(R)^G$ as a subring, where $\mathcal{Z}(R)^G=\{r \in \mathcal{Z}(R)~|~ \sigma_{g}(r)=r~\forall~ g \in G\}$ is the fixed subring of $\mathcal{Z}(R)$. The following result provides a neccesary condition for the equality $Z(R*^{\sigma}_{\tau}G)= \mathcal{Z}(R)^G$ to hold:
\begin{lemma}{\label{lemma1}}{\rm{(\cite{Jes}, Lemma 2.6.1)}}
Let $R*_{\tau}^{\sigma}G$ be a crossed product with action $\sigma$ and assume that $R$ is simple. If $\sigma_{g}$ is not inner in $R$ for every non-trivial $g \in G,$ then $R*_{\tau}^{\sigma}G$ is a simple ring and $$\mathcal{Z}(R*_{\tau}^{\sigma}G)=\mathcal{Z}(R)^{G}.$$
\end{lemma}
\noindent A \textit{classical crossed product} is a crossed product $L*_{\tau}^{\sigma} G$, where $L$ is a finite Galois extension of $F=\mathcal{Z}(L*_{\tau}^{\sigma}G)$, $G=\G(L/F)$ is the Galois group of $L$ over $F$ and $\sigma$ is the natural action of $G$ over $L$. This classical crossed product $L*_{\tau}^{\sigma} G$ is commonly denoted by $(L/F,\tau).$ A \textit{cyclotomic algebra} over $F$ is a classical crossed product $(L/F,\tau),$ where $L=F(\zeta)$ is a cyclotomic extension over $F$. Here $\zeta$ is a root of unity. A \textit{cyclic algebra} is a classical crossed product $(L/F,\tau)$, where $L$ over $F$ is a cyclic extension, i.e, $\G(L/F)$ is cyclic. If $\G(L/F)$ is generated by $\phi$ and has order $d$, then $(L/F,\tau)$ is isomorphic to $(L/F,\tau')$, where $\tau'$ is given by $$ \tau'(\phi^i,\phi^j) = \left\{ 
	                                 \begin{array}{ll}
	                                 1, & i+j <d;\\
	                                 z{_\phi}^{d}, &  i+j \geq d.
	                                 \end{array}
\right.
$$ Such a cyclic algebra $(L/F,\tau)$ is commonly denoted by $(L/F, \phi, z_{\phi}^d).$ For more details on crossed products, see (\cite{Jes}, Section 2.6).
\para \noindent The following is a classical result (\cite{Jes}, Theorem 3.7.1) due to Brauer and Witt:
\begin{theorem}\label{thm1} {\rm(Brauer-Witt theorem)}
	Every Wedderburn component of a semisimple group algebra is Brauer equivalent to a cyclotomic algebra over its center.
\end{theorem}\noindent We are now going to recall Proposition 3.4 of \cite{OdRS04} by Olivieri, del R{\'{\i}}o and Sim{\'o}n, which we wish to generalize. For a Shoda pair $(H,K)$ of $G$, let  $$\widehat{H}:=\frac{1}{|H|}\displaystyle\sum_{h \in H}h,$$ $$\varepsilon(H,K):=\left\{\begin{array}{ll}\widehat{K}, & \hbox{$H=K$;} \\\prod(\widehat{K}-\widehat{L}), & \hbox{otherwise,}\end{array}\right.$$ where $L$ runs over the normal subgroups of $H$ minimal with respect to the property of including $K$ properly, and $$e(G,H,K):= {\rm~the~sum~of~all~the~distinct~}G{\rm {\operatorname{-}} conjugates~of~}\varepsilon(H,K).$$ \para \noindent
A Shoda pair $(H,K)$ of $G$ is called a \textit{strong Shoda pair} of $G$ (see \cite{OdRS04}, Definition 3.1, Proposition 3.3) if the following conditions hold:\begin{description}\item [(i)]$H$ is normal in $\operatorname{Cen}_{G}(\varepsilon(H,K))$, the centralizer of $\varepsilon(H,K)$ in $G$;\item [(ii)] the distinct $G$-conjugates of $\varepsilon(H,K)$ are mutually orthogonal.\end{description} In \cite{OdRS04}, it is proved that if $(H,K)$ is a strong Shoda pair of $G$ and $\lambda$ a linear character on $H$ with kernel $K$, then $e_{\Q}(\lambda^G)= e(G,H,K).$ \begin{theorem}\label{thm2} {\rm(\cite{OdRS04}, Proposition 3.4)} Let $(H,K)$ be a strong Shoda pair of $G$ and let $k=[H:K]$. Let $N=N_{G}(K)$, the normalizer of $K$ in $G$ and $n=[G:N]$. Let $x$ be a generator of $H/K$ and $\phi: N/H \rightarrow N/K$ a left inverse of the projection $N/K \rightarrow N/H$. Then 
	$$ \Q Ge(G,H,K)\cong M_{n}(\Q(\zeta_{k})*^{\sigma}_{\tau}N/H),$$
	where $\zeta_{k}$ is a primitive $k$-th root of unity, the action $\sigma$ and the twisting $\tau$ are given by $\sigma_{a}(\zeta_{k})=\zeta_{k}^{i}$ if ${\phi(a)}^{-1}x \phi(a)=x^{i};$ $\tau(a,b)=\zeta_{k}^{j}$ if $\phi(ab)^{-1}\phi(a)\phi(b)=x^{j}$, for $a,b \in N/H$ and integers $i$ and $j$.
\end{theorem}
\noindent A Shoda  pair $(H,K)$  of $G$ is  called a \textit{generalized strong Shoda pair} of $G$ (\cite{BK3}, Section 2) if there is a chain $H=H_{0}\leq H_{1}\leq \cdots \leq H_{n}=G$ of subgroups of $G$ such that the following conditions hold for all $ 0 \leq i \leq n-1$: \begin{enumerate}[(a)] \item $H_i \unlhd \operatorname{Cen}_{H_{i+1}}(e_{\mathbb{Q}}(\lambda^{H_i}))$; \item the distinct $H_{i+1}$-conjugates of $e_{\mathbb{Q}}(\lambda^{H_i})$ are mutually orthogonal.\end{enumerate} Here $\lambda$  is a linear character on $H$ with kernel $K$. \para \noindent For a generalized strong Shoda pair $(H,K)$ of $G$, a chain $H=H_{0}\leq H_{1}\leq \cdots \leq H_{n}=G$ of subgroups of $G$ which satisfy  the conditions (a) and (b) above, for all $0 \leq i \leq n-1$,  is  called a \textit{strong inductive chain} of the generalized strong Shoda pair $(H,K)$. \para \noindent A monomial character $\chi$ of $G$ is called {\it generalized strongly monomial} if it arises from a generalized strong Shoda pair $(H,K)$.
\begin{remark}\label{remark2}
(i) If $(H,K) $ is a Shoda pair and $H=H_{0}\leq H_{1}\leq \cdots \leq H_{n}=G$ are subgroups of $G$ such that $H_j \unlhd H_{j+1} $ for some $j$, then the conditions (a) and (b) stated above hold true for $i=j$. The condition (a) is obvious and (b) holds true because for any $g \in H_{j+1}\backslash \operatorname{Cen}_{H_{j+1}}(e_{\Q}(\lambda^{H_j}))$, $e_{\Q}(\lambda^{H_j})$ and $(e_{\Q}(\lambda^{H_j}))^g$ being distinct primitive central idempotents of $\Q H_{j}$ are mutually orthogonal.\\
(ii) From part (i), it also follows that if $(H,K)$ is a Shoda pair with $H \unlhd G,$ then $(H,K)$ is a  strong Shoda pair.
\end{remark} \noindent  When $(H,K)$ is a generalized strong Shoda pair of $G$ and $H=H_{0}\leq H_{1}\leq \cdots \leq H_{n}=G$ is a strong inductive chain from $H$ to $G$, then it is  proved in (\cite{BK3}, Proposition 2) that, for all $0 \leq i \leq n-1$, \begin{equation}{\label{eqn}}
\Q H_{i+1}e_{\mathbb{Q}}(\lambda^{H_{i+1}})  \cong M_{k_i}(\Q H_{i}e_{\mathbb{Q}}(\lambda^{H_{i}}) *^{\sigma _{H_{i}}}_{\tau_{H_{i}}} C_i/H_i),
\end{equation}  where $k_{i}=[H_{i+1}:C_{i}], ~C_{i}=\operatorname{Cen}_{H_{i+1}}(e_{\Q}(\lambda^{H_i})), ~ \sigma_{H_i} : C_{i}/H_{i} \rightarrow \operatorname{Aut}(\Q H_{i}e_{\Q}(\lambda^{H_{i}}))$ maps $x_{i} \in C_{i}/H_{i}$ to the conjugation automorphism $(\sigma_{H_{i}})_{x_{i}}$ on $\Q H_{i}e_{\Q}(\lambda^{H_{i}})$ induced by $\overline{x}_{i}$, i.e., $(\sigma_{H_{i}})_{x_{i}}(\alpha)= \overline{x}_{i}^{-1}\alpha \overline{x}_{i}$. Here $\overline{x}_{i} \in C_{i}$ is a fixed representative of the coset $x_{i}$. Also $\tau_{H_{i}}: C_{i}/H_{i} \times C_{i}/H_{i} \rightarrow \mathcal{U}(\Q H_{i}e_{\Q}(\lambda^{H_{i}}))$ is given by $\tau_{H_{i}}(x_{i},y_{i}) = \overline{x_{i}y_{i}}^{-1}\overline{x}_{i}\overline{y}_{i}e_{\Q}(\lambda^{H_{i}}),$ for $x_{i},y_{i} \in C_{i}/H_{i}.$ Furthermore, the isomorphism in equation (\ref{eqn}) maps $\alpha \in \Q H_{i+1}e_{\Q}(\lambda^{H_{i+1}})$ to $(\alpha_{ml})_{k_{i} \times k_{i}}$, where $\alpha_{ml}= e_{\Q}(\lambda^{H_{i}})t_{l}\alpha t_{m}^{-1}e_{\Q}(\lambda^{H_{i}})$ and $\{t_{i},t_{2},\cdots,t_{k_{i}}\}$ is a transversal set of $C_{i}$ in $H_{i+1}.$
\para \noindent The following structure of the simple component of $\mathbb{Q}G$ realized by a generalized strong Shoda pair is  proved in (\cite{BK3}, Theorem 3):
\begin{theorem}\label{thm3} Let $(H,K)$ be a generalized strong Shoda pair of $G$ and $\lambda$ a linear character on $H$ with kernel $K$.  Let $H=H_{0}\leq H_{1}\leq \cdots \leq H_{n}=G$ be a strong inductive chain from $H$ to $G$. Then {\footnotesize $$\mathbb{Q}Ge_{\mathbb{Q}}(\lambda^{G}) \cong  {M_{k_{n-1}}(\cdots(M_{k_{0}}(\Q H e_{\Q}(\lambda)
			*^{\sigma_{H_{0}}}_{\tau_{H_{0}}} C_{0}/H_{0}) *^{\sigma_{H_{1}}}_{\tau_{H_{1}}} \cdots) *^{\sigma_{H_{n-1}}}_{\tau_{H_{n-1}}} C_{n-1}/H_{n-1})},$$} where $C_{i},~\sigma_{H_{i}},~\tau_{H_{i}},~k_{i}$ are as defined above. \end{theorem}
\section{Main Theorem}
We begin by fixing some notations to be used throughout the paper. Given a generalized strong Shoda pair $(H,K)$ of $G,$   by $H=H_{0}\leq H_{1} \leq \cdots \leq H_{n} = G$, we always mean a strong inductive chain from $H$ to $G,$ and further set the following:
\begin{equation*}
\begin{array}{lll} \noindent
\lambda&:=& \mbox{a linear character on H with kernel K.}\\ 
C_{i}& :=& \C_{H_{i+1}}(e_{\Q}(\lambda^{H_{i}})), 0 \leq i \leq n-1.\\
k_{i} &:=& \mbox{the index of}~ C_{i}~ \mbox{in}~ H_{i+1} , 0 \leq i \leq n-1.\\
k &:=& k_{0}k_{1} \cdots k_{n-1}.\\ 
(\sigma_{H_{i}})_{x_{i}} &:=& \mbox{the automorphism of}~ \Q H_{i}e_{\Q}(\lambda^{H_{i}})~ \mbox{which sends}~\alpha ~\mbox{to}~  \overline{x}_{i}^{-1}\alpha\overline{x}_{i},~\mbox{where} \\
&&\alpha \in \Q H_{i}e_{\Q}(\lambda^{H_{i}})~\mbox{and} ~\overline{x}_{i} \in C_{i} ~\mbox{is a fixed representative of the coset}~ \\&& x_{i} \in C_{i}/H_{i},0 \leq i \leq n-1.\\
(\tilde{\sigma}_{H_{i}})_{x_{i}} &:=& \mbox{the restriction of}~ (\sigma_{H_{i}})_{x_{i}}~\mbox{on the center}~\mathcal{Z}(\Q H_{i}e_{\Q}(\lambda^{H_{i}}))~\mbox{of}~\Q H_{i}e_{\Q}(\lambda^{H_{i}}).\\
\mathcal{A} &:=& M_{k_{n-1}}(\cdots(M_{k_{0}}(\Q He_{Q}(\lambda) *^{\sigma_{H_{0}}}_{\tau_{H_{0}}} C_{0}/H_{0}) *^{\sigma_{H_{1}}}_{\tau_{H_{1}}} \cdots)*^{\sigma_{H_{n-1}}}_{\tau_{H_{n-1}}} C_{n-1}/H_{n-1}),
\mbox{the} \linebreak\\&& \mbox{algebra which appeared in Theorem \ref{thm3}}.\\
\mathbb{F} &:=& \mathcal{Z}(\mathcal{A}).\\
\mbox{E} &:=& \Q He_{\Q}(\lambda).\\
\mathbb{E} &:=& k \times k~ \mbox{scalar matrices}~ \operatorname{diag}(\alpha,\alpha,\cdots,\alpha)_{k}, \alpha \in \mbox{E}.\\
\mbox{F} &:=&\{\alpha \in \mbox{E} ~|~ \operatorname{diag}(\alpha,\alpha,\cdots,\alpha)_{k} \in \mathbb{F}\}.\\
\mathcal{B} &:=& M_{k \times k}(\mbox{F}).
\end{array}
\end{equation*}

\noindent
In section 2, we pointed out in Theorem \ref{thm3} that the simple component $\Q Ge_{\Q}(\lambda^G)$ of $\Q G$ is isomorphic to $\mathcal{A}.$ In this section, we will explain the structure of $\mathcal{A}.$ The following steps shall be proved:  
\begin{equation*}
\begin{array}{lll}
\mbox{Step (i)} &:& \mathcal{A} \cong M_{k \times k}(\operatorname{Cen}_{\mathcal{A}}(\mathcal{B})).\\
\mbox{Step (ii)} &:& \mathbb{E}~\mbox{is a Galois extension of} ~\mathbb{F}~ \mbox{and the}~ \operatorname{dim}_{\mathbb{F}}(\mathbb{E})= \prod_{i=0}^{n-1}|C_{i}/H_{i}|.\\ 
\mbox{Step (iii)} &:& \operatorname{Cen}_{\mathcal{A}}(\mathcal{B})~\mbox{contains}~ \mathbb{E}~ \mbox{as a maximal subfield}.\\
\mbox{Step (iv)} &:& \mbox{Finally,}~ \operatorname{Cen}_{\mathcal{A}}(\mathcal{B})~\mbox{is isomorphic to the classical crossed product}\\&&
(\mathbb{E}/\mathbb{F},\tau)~\mbox{with respect to a factor set} ~\tau.
\end{array}
\end{equation*}
\noindent
\textbf{Proof of Step (i).} The center of $\mathcal{B}$ is clearly $\mathbb{F}.$ Hence, $\mathcal{B}$ is a central $\mathbb{F}$- subalgebra of the central simple $\mathbb{F}$-algebra $\mathcal{A}.$ By applying Double Centralizer Theorem (\cite{Jes}, Theorem 2.1.10), we obtain that 
$$\mathcal{A} \cong \mathcal{B} \otimes_{\mathbb{F}} \operatorname{Cen}_{\mathcal{A}}(\mathcal{B}).$$ 
It can be checked that the map $\mathcal{B} \otimes \operatorname{Cen}_{\mathcal{A}}(\mathcal{B}) \rightarrow M_{k \times k}(\operatorname{Cen}_{\mathcal{A}}(\mathcal{B}))$ given by $(\alpha_{ij}) \otimes x \mapsto (\beta_{ij}),$ where $(\alpha_{ij})_{k \times k} \in \mathcal{B}, x \in \operatorname{Cen}_{\mathcal{A}}(\mathcal{B})$ and $\beta_{ij} = \operatorname{diag}(\alpha_{ij},\alpha_{ij},\cdots,\alpha_{ij})_{k}x$ is an isomorphism. This proves step (i).\qed \para \noindent
\textbf{Proof of Step (ii).} For notational convenience, let us denote the central simple algebra $\small{M_{k_{n-2}}(\cdots(M_{k_{0}}(\Q H e_{\Q}(\lambda)*^{\sigma_{H_{0}}}_{\tau_{H_{0}}} C_{0}/H_{0}) *^{\sigma_{H_{1}}}_{\tau_{H_{1}}} \cdots) *^{\sigma_{H_{n-2}}}_{\tau_{H_{n-2}}} C_{n-2}/H_{n-2})}$ by $\mathcal{A}_{1}.$ We have
\begin{equation*}
\mathcal{A} = M_{k_{n-1}}(\mathcal{A}_{1} *^{\sigma_{H_{n-1}}}_{\tau_{H_{n-1}}} C_{n-1}/H_{n-1}).
\end{equation*} 
By Theorem \ref{thm3}, $\mathcal{A}_{1}$ is isomorphic to the simple component $\Q H_{n-1} e_{\Q}(\lambda^{H_{n-1}})$ of the group algebra $\Q H_{n-1}.$ Hence, $\mathcal{A}_{1}$ is a simple ring. Also, by Theorem (\cite{BK5}, Proposition 1 (a)), $(\sigma_{H_{n-1}})_{x}$ is not inner for any $1 \neq x \in C_{n-1}/H_{n-1}.$ Therefore, Lemma \ref{lemma1} implies that $$\mathcal{Z}(\mathcal{A}_{1} *^{\sigma_{H_{n-1}}}_{\tau_{H_{n-1}}} C_{n-1}/H_{n-1})= \mathcal{Z}(\mathcal{A}_{1})^{C_{n-1}/H_{n-1}}$$ and hence,
\begin{equation}{\label{eq1}}
\mathcal{Z}(\mathcal{A}) =\{ \operatorname{diag}(\alpha,\alpha,\cdots,\alpha)_{k_{n-1}} |~ \alpha \in \mathcal{Z}(\mathcal{A}_{1})^{C_{n-1}/H_{n-1}} \}.
\end{equation}
Now, if $\mathcal{A}_{2} = \small{M_{k_{n-3}}(\cdots(M_{k_{0}}(\Q H e_{\Q}(\lambda)*^{\sigma_{H_{0}}}_{\tau_{H_{0}}} C_{0}/H_{0}) *^{\sigma_{H_{1}}}_{\tau_{H_{1}}} \cdots)*^{\sigma_{H_{n-3}}}_{\tau_{H_{n-3}}} C_{n-3}/H_{n-3})},$ then 
$\mathcal{A}_{1} = M_{k_{n-2}}(\mathcal{A}_{2} *^{\sigma_{H_{n-2}}}_{\tau_{H_{n-2}}} C_{n-2}/H_{n-2})$
and 
\begin{equation}{\label{eq2}}
\mathcal{Z}(\mathcal{A}_{1}) = \{ \operatorname{diag}(\alpha, \alpha, \cdots , \alpha)_{k_{n-2}} | ~ \alpha \in \mathcal{Z}(\mathcal{A}_{2})^{C_{n-2}/H_{n-2}} \}.
\end{equation}
Equations (\ref{eq1}) and (\ref{eq2}) together give us 
$$ \mathcal{Z}(\mathcal{A}) = \{ \operatorname{diag}(\alpha,\alpha,\cdots,\alpha)_{k_{n-1}k_{n-2}} ~|~ \alpha \in (\mathcal{Z}(\mathcal{A}_{2})^{C_{n-2}/H_{n-2}})^{C_{n-1}/H_{n-1}} \}.$$
Applying the same process to $\mathcal{A}_{2}$ and continue, we finally obtain that $$\mathbb{F} = \mathcal{Z}(\mathcal{A}) = \{ \operatorname{diag}(\alpha,\alpha,\cdots,\alpha)_{k} ~|~ \alpha \in ((\Q He_{\Q}(\lambda)^{C_{0}/H_{0}}) \cdots)^{C_{n-1}/H_{n-1}} \}$$ and hence it is contained in $\mathbb{E}.$ Now as $\mathbb{E} \cong \Q(\zeta_{n}),$ where $n=[H:K]$, is a Galois extension of $\Q$ and $\Q \subseteq \mathbb{F} \subseteq \mathbb{E}$, it turns out that $\mathbb{E}$ is a Galois extension over $\mathbb{F}$. \para \noindent We will now show that $\operatorname{dim}_{\mathbb{F}}(\mathbb{E})=\prod_{i=0}^{n-1}|C_{i}/H_{i}|.$ Let $\mathbb{E}_{0}= \mathbb{E},~ \mathbb{E}_{n}= \mathbb{F}$ and for $0 < i < n,$ let \begin{equation}{\label{eq}}
\mathbb{E}_{i}= \{ \operatorname{diag}(\alpha,\alpha.\cdots,\alpha)_{k_{i}k_{i+1}\cdots k_{n-1}} ~|~ \alpha \in \mathcal{Z}(\mathcal{A}_{n-i})\},
\end{equation} where $\mathcal{A}_{n-i} = M_{k_{i-1}}(\cdots(M_{k_{0}}(\Q He_{\Q}(\lambda) *^{\sigma_{H_{0}}}_{\tau_{H_{0}}} C_{0}/H_{0})*^{\sigma_{H_{1}}}_{\tau_{H_{1}}} \cdots) *^{\sigma_{H_{i-1}}}_{\tau_{H_{i-1}}} C_{i-1}/H_{i-1})).$ \para \noindent
\textbf{Claim}: $\mathbb{E} = \mathbb{E}_{0} \supseteq \mathbb{E}_{1} \supseteq \cdots \supseteq \mathbb{E}_{n}= \mathbb{F}$ are subfields of the extension $\mathbb{E}/\mathbb{F}$ and the $\operatorname{Gal}(\mathbb{E}_{i}/\mathbb{E}_{i+1})=\{ (\tilde{\sigma}_{H_{i}})_{x_{i}} ~|~ x_{i} \in C_{i}/H_{i} \}$ is isomorphic to $C_{i}/H_{i},~ 0 \leq i \leq n-1.$ \para \noindent
From equation (\ref{eq1}), $\mathbb{E}_{n}= \mathcal{Z}(\mathcal{A})= \{ \operatorname{diag}(\alpha,\alpha,\cdots,\alpha)_{k_{n-1}}~|~\alpha \in \mathcal{Z}(\mathcal{A}_{1})^{C_{n-1}/H_{n-1}} \}.$ The right hand side of above equation is clearly contained in $\mathbb{E}_{n-1}.$ Therefore, $\mathbb{E}_{n} \subseteq \mathbb{E}_{n-1}.$ Similarly, $\mathbb{E}_{n-1} \subseteq \mathbb{E}_{n-2} \subseteq \cdots \subseteq \mathbb{E}_{0}=\mathbb{E}.$ For $0 \leq i \leq n-1,$ it is shown in (\cite{BK5}, Proposition 1, Corollary 1) that the map $C_{i}/H_{i} \rightarrow \operatorname{Aut}(\mathbb{E}_{i})$ given by $x \mapsto (\tilde{\sigma}_{H_{i}})_{x}$ is a one-one group homomorphism and $\mathbb{E}_{i+1}$ is the fixed subfield of $\mathbb{E}_{i}$ under the action of $C_{i}/H_{i}.$ This gives that $C_{i}/H_{i} \cong \G(\mathbb{E}_{i}/\mathbb{E}_{i+1})$ for all $0 \leq i \leq n-1$ and hence proves the claim which immediately yields that $\operatorname{dim}_{\mathbb{F}}(\mathbb{E})=\prod_{i=0}^{n-1}|C_{i}/H_{i}|$ follows. \qed \para \noindent
\textbf{Proof of Step (iii).} Recall that $$\mathcal{A}={M_{k_{n-1}}(\cdots(M_{k_{0}}(\Q H e_{\Q}(\lambda)*^{\sigma_{H_{0}}}_{\tau_{H_{0}}} C_{0}/H_{0}) *^{\sigma_{H_{1}}}_{\tau_{H_{1}}} \cdots) *^{\sigma_{H_{n-1}}}_{\tau_{H_{n-1}}} C_{n-1}/H_{n-1})}$$ and clearly $\operatorname{dim}_{\mathbb{F}}(\mathcal{A}) = k_{n-1}^2 \cdots k_{0}^2 |C_{n-1}/H_{n-1}| \cdots |C_{0}/H_{0}| \operatorname{dim}_{\mathbb{F}}(\mathbb{E}).$ Also note that $\C_{\mathcal{A}}(\mathcal{B})$ contains $\mathbb{E}.$ By step (ii), $\operatorname{dim}_{\mathbb{F}}(\mathbb{E}) = |C_{0}/H_{0}| |C_{1}/H_{1}| \cdots |C_{n-1}/H_{n-1}|.$ Hence, 
\begin{equation*}
\operatorname{dim}_{\mathbb{F}}(\mathcal{A}) = k^2 |C_{0}/H_{0}|^2 \cdots |C_{n-1}/H_{n-1}|^2.
\end{equation*} 
However by step (i), 
\begin{equation*}
\operatorname{dim}_{\mathbb{F}}(\mathcal{A}) = k^2 \operatorname{dim}_{\mathbb{F}}(\C_{\mathcal{A}}(\mathcal{B})).
\end{equation*}
Consequently,  
$$ \operatorname{dim}_{\mathbb{F}}(\C_{\mathcal{A}}(\mathcal{B})) = |C_{0}/H_{0}|^2 |C_{1}/H_{1}|^2 \cdots |C_{n-1}/H_{n-1}|^2 = (\operatorname{dim}_{\mathbb{F}} \mathbb{E})^2,$$ which, by (\cite{Jes}, Proposition 2.2.4 (2)), implies that $\mathbb{E}$ is maximal subfield of $\C_{\mathcal{A}}(\mathcal{B})).$ \qed \para \noindent
Step (iv) follows naturally from step (iii), which is obvious to the experts. However, for the reader's convenience, we outline its proof below:\para \noindent
\textbf{Proof of Step (iv).} By Noether-Skolem theorem (\cite{Jes}, Theorem 2.1.9), for every $\sigma \in \G(\mathbb{E}/\mathbb{F}),$ there exist an invertible $z_{\sigma}$ in $\C_{\mathcal{A}}(\mathcal{B})$ such that $\sigma(\alpha)= z_{\sigma}^{-1} \alpha z_{\sigma}$ for all $\alpha \in \mathbb{E}.$ It turns out that $\{z_{\sigma} ~|~ \sigma \in \G(\mathbb{E}/\mathbb{F})\}$ are linearly independent over $\mathbb{E}.$ As $\operatorname{dim}_{\mathbb{E}}(\C_{\mathcal{A}}(\mathcal{B}))= |\G(\mathbb{E}/\mathbb{F})|,$ it follows that $\{z_{\sigma} ~|~ \sigma \in \G(\mathbb{E}/\mathbb{F})\}$ is a basis for the $\mathbb{E}$-vector space $\C_{\mathcal{A}}(\mathcal{B})$. If $\tau: \G(\mathbb{E}/\mathbb{F}) \times \G(\mathbb{E}/\mathbb{F}) \rightarrow \mathcal{U}(\mathbb{E})$ is given by $(\sigma,\tau) \mapsto z_{\sigma\tau}^{-1} z_{\sigma} z_{\tau},$ then the associativity of multiplication in $\operatorname{Cen}_{\mathcal{A}}(\mathcal{B})$ implies that $\tau$ is a factor set and $\C_{\mathcal{A}}(\mathcal{B})$ is isomorphic to the classical crossed product $(\mathbb{E}/\mathbb{F},\tau).$ \qed \para \noindent
Steps (i)-(iv) yield the following theorem:
\begin{theorem}{\label{thm4}}
	Let $(H,K)$ be a generalized strong Shoda pair of $G$ and $\lambda$ a linear character on $H$ with kernel $K.$ Then the simple component $\Q Ge_{\Q}(\lambda^G)$ of $\Q G$ is isomorphic to a matrix algebra over the classical crossed product $(\mathbb{E}/\mathbb{F},\tau)$ for some twisting $\tau,$ where $\mathbb{E} \cong \Q He_{\Q}(\lambda) \cong \Q(\lambda)$ and $\mathbb{F} \cong \Q(\lambda^G).$ 
\end{theorem} \noindent
In general for a given irreducible character $\chi$ of $G$, it is well known that $m_{\Q}(\chi)$ divides $\chi(1)$. The following corollary is an improvement of this result when $\chi$ is a generalized strongly monomial character and follows immediately from the above theorem.
\begin{cor}
	If $(H,K)$ is a generalized strong Shoda pair and $\lambda$ a linear character on $H$ with kernel $K$, then $m_{\Q}(\lambda^G)$ divides $\prod_{i=0}^{n-1}|C_{i}/H_{i}|$, where $C_{i}$'s and $H_{i}$'s are as described above.
\end{cor}
\subsection{Understanding $\G(\mathbb{E}/\mathbb{F})$ and $\tau$}
In order to understand the Schur index of the algebra $\Q Ge_{\Q}(\lambda^G),$ one further needs to understand Galois group $\G(\mathbb{E}/\mathbb{F})$ and the twisting $\tau.$ This can be made clear by determining a set of basis units $\{z_{\sigma} ~|~ \sigma \in \G(\mathbb{E}/\mathbb{F})\}$ of $(\mathbb{E}/\mathbb{F},\tau)$ which appeared in the proof of step (iv). \para \noindent
A group $G$ is said to be of type $G_{1}$-by-$G_{2}$ if $G$ has a normal subgroup isomorphic to $G_{1}$ and the quotient by $G_{1}$ is isomorphic to $G_{2}.$ From the proof of step (ii), we understand that $\G(\mathbb{E}/\mathbb{F})$ is of the type $(((C_{0}/H_{0}$-by-$C_{1}/H_{1})$-by-$C_{2}/H_{2})\cdots)$-by-$C_{n-1}/H_{n-1}.$ For $x_{i} \in C_{i}/H_{i},$ if we consider an arbitrary lifting, say $\sigma_{x_{i}},$ of $(\tilde{\sigma}_{H_{i}})_{x_{i}}$ to $\mathbb{E},$ then it turns out that $$\G(\mathbb{E}/\mathbb{F})= \{\sigma_{x_{0}}\circ \sigma_{x_{1}} \circ \cdots \circ \sigma_{x_{n-1}} ~|~ x_{i} \in C_{i}/H_{i}, 0 \leq i \leq n-1\}.$$ For any $0 \leq i \leq n-1$ and $x_{i} \in C_{i}/H_{i}$, by Noether-Skolem Theorem (\cite{Jes}, Theorem 2.1.9), we can find some $z_{x_{i}}$ in $\operatorname{Cen}_{\mathcal{A}}(\mathcal{B}),$ so that $\sigma_{x_{i}}$ is given by conjugation on $\mathbb{E}$ by $z_{x_{i}},$ i.e., $\sigma_{x_{i}}(\alpha)= z_{x_{i}}^{-1} \alpha z_{x_{i}}~\forall~\alpha \in \mathbb{E}.$ It can be checked that if $\sigma = \sigma_{x_{0}} \circ \sigma_{x_{1}} \circ \cdots \sigma_{x_{n-1}},$ then by setting $z_{\sigma}= z_{x_{0}} z_{x_{1}} \cdots z_{x_{n-1}},$ it follows that $\sigma(\alpha)= z_{\sigma}^{-1} \alpha z_{\sigma} ~ \forall ~\alpha \in \mathbb{E}.$ Hence, $\{ z_{x_{0}}z_{x_{1}}\cdots z_{x_{n-1}}~|~ x_{i} \in C_{i}/H_{i}, 0 \leq i \leq n-1 \}$ turns out to be a set of basis units of $\operatorname{Cen}_{\mathcal{A}}(\mathcal{B})$ and the twisting $\tau:  \G(\mathbb{E}/\mathbb{F}) \times \G(\mathbb{E}/\mathbb{F}) \rightarrow \mathcal{U}(\mathbb{E}),$ which appeared in Theorem \ref{thm4} is given as follows:
$$\tau(\sigma_{x_{0}} \circ \cdots \circ \sigma_{x_{n-1}},\sigma_{y_{0}} \circ \cdots \circ \sigma_{y_{n-1}}) = (z_{x_{0}y_{0}} \cdots  z_{x_{n-1}y_{n-1}})^{-1}z_{x_{0}} \cdots z_{x_{n-1}} z_{y_{0}} \cdots z_{y_{n-1}}.$$ 
The problem thus boils down to determine $z_{x_{i}},$ for $x_{i} \in C_{i}/H_{i},$ which we wish to do next.\para \noindent
\textbf{(3.1.1)} {\it If} $z \in \operatorname{Cen}_{\mathcal{A}}(\mathcal{B})$ {\it is such that}  $z^{-1}\mathbb{E}z=\mathbb{E},$ {\it then the conjugation on} $\mathbb{E}$ {\it by} $z$ {\it extends} $(\tilde{\sigma}_{H_{i}})_{x_{i}},$ {\it for some} $x_{i} \in C_{i}/H_{i},~0 \leq i \leq n-1,$ {\it if and only if} $z =  \operatorname{diag}(\overline{x}_{i} u_{x_{i}}, \overline{x}_{i} u_{x_{i}},
\cdots, \overline{x}_{i} u_{x_{i}})_{k_{i} \cdots k_{n-1}}$ {\it for some} $u_{x_{i}}$ {\it in} $\mathcal{A}_{n-i}.$ \para
\noindent
Let $\sigma \in \G(\mathbb{E}/\mathbb{F})$ be given by $\alpha \mapsto z^{-1}\alpha z,~\alpha \in \mathbb{E}.$ We are given that $\sigma$ extends $(\tilde{\sigma}_{H_{i}})_{x_{i}}$ for some $x_{i} \in C_{i}/H_{i}.$ This means that for any $\alpha \in \mathbb{E}_{i},$
$$z^{-1} \alpha z = X_{i}^{-1} \alpha X_{i},$$ where $X_{i}=\operatorname{diag}(\overline{x}_{i},\overline{x}_{i},\cdots,\overline{x}_{i})_{k_{i}\cdots k_{n-1}}.$ Therefore, $zX_{i}^{-1} \in \C_{\mathcal{A}}(\mathbb{E}_{i}).$ Let us see that $\operatorname{Cen}_{\mathcal{A}}(\mathbb{E}_{i}) = M_{k_{i}\cdots k_{n-1}}(\mathcal{A}_{n-i}).$ Clearly, 
\begin{equation}{\label{eq3}}
M_{k_{i}\cdots k_{n-1}}(\mathcal{A}_{n-i}) \subseteq \operatorname{Cen}_{\mathcal{A}}(\mathbb{E}_{i}).
\end{equation}
By (\cite{Jes}, Proposition 2.2.4 (1)), $$\operatorname{dim}_{\mathbb{E}_{i}}(\operatorname{Cen}_{\mathcal{A}}(\mathbb{E}_{i})) = \frac{\operatorname{dim}_{\mathbb{F}}(\mathcal{A})}{[\mathbb{E}_{i}:\mathbb{F}]^2}= \frac{k_{0}^2k_{1}^2\cdots k_{n-1}^2 |C_{0}/H_{0}|^2 |C_{1}/H_{1}|^2 \cdots |C_{n-1}/H_{n-1}|^2}{|C_{i}/H_{i}|^2 \cdots |C_{n-1}/H_{n-1}|^2}.$$  Also $\operatorname{dim}_{\mathbb{E}_{i}}(M_{k_{i}\cdots k_{n-1}}(\mathcal{A}_{n-i})) = k_{n-1}^2 \cdots k_{0}^2 |C_{0}/H_{0}|^2 \cdots |C_{i-1}/H_{i-1}|^2,$ which gives equality in equation (\ref{eq3}). Hence, we can write 
$$zX_{i}^{-1} = (A_{mn})_{l \times l},$$ where $l = k_{i}k_{i+1} \cdots k_{n-1}$ and $A_{mn} \in \mathcal{A}_{n-i}$ for $1 \leq m,n \leq l.$ \para \noindent Let $e$ be $\operatorname{diag}(e_{\Q}(\lambda),e_{\Q}(\lambda),\cdots,e_{\Q}(\lambda))_{k_{0}k_{1}\cdots k_{i-1}}$ and $E_{ij}$ be $l \times l$ block matrix whose $(i,j)$th block is $e$ and all other blocks are zero. As $z \in \operatorname{Cen}_{\mathcal{A}}(\mathcal{B})$ and $E_{ij} \in \mathcal{B},$ we have $zE_{ij}=E_{ij}z,$ which gives 
$$(A_{mn}\overline{x}_{i}) E_{ij} = E_{ij} (A_{mn}\overline{x}_{i}).$$
As the above equation is true for all $E_{ij},~1 \leq i,j \leq l,$ it yields that, $A_{mn}\overline{x}_{i}=0$ if $m \neq n$ and $A_{mm}\overline{x}_{i}=A_{nn}\overline{x}_{i}$ for all $1 \leq m,n \leq l.$ If we denote $\overline{x}_{i}^{-1}A_{11}\overline{x}_{i}$ by $u_{x_{i}},$ then it turns out that $z = \operatorname{diag}(\overline{x}_{i}u_{x_{i}},\overline{x}_{i}u_{x_{i}},\cdots,\overline{x}_{i}u_{x_{i}})_{k_{i} \cdots k_{n-1}}.$ The invertibility of $u_{x_{i}}$ follows from that of $z$ and $\overline{x}_{i}.$
\para \noindent Conversely, if $z = \operatorname{diag}(\overline{x}_{i}u_{x_{i}},\overline{x}_{i}u_{x_{i}},\cdots,\overline{x}_{i}u_{x_{i}})_{k_{i} \cdots k_{n-1}}$ and $X = \operatorname{diag}(\alpha,\alpha,\cdots,\alpha)_{k_{i} \cdots k_{n-1}} \in \mathbb{E}_{i},$ where $\alpha \in \mathcal{Z}(\mathcal{A}_{n-i}),$ then 
\begin{flalign*}
\sigma(X)= z^{-1}Xz &= (\operatorname{diag}(\overline{x}_{i}u_{x_{i}},\overline{x}_{i}u_{x_{i}},\cdots,\overline{x}_{i}u_{x_{i}}))^{-1} X \operatorname{diag}(\overline{x}_{i}u_{x_{i}},\overline{x}_{i}u_{x_{i}},\cdots,\overline{x}_{i}u_{x_{i}})\\&
= \operatorname{diag}(u_{x_{i}}^{-1}\overline{x}_{i}^{-1}\alpha \overline{x}_{i}u_{x_{i}}, \cdots, u_{x_{i}}^{-1}\overline{x}_{i}^{-1}\alpha \overline{x}_{i}u_{x_{i}})\\&
= \operatorname{diag}(\overline{x}_{i}^{-1}\alpha \overline{x}_{i},\cdots,\overline{x}_{i}^{-1}\alpha \overline{x}_{i}).
\end{flalign*}
The second last equality follows because $\overline{x}_{i}^{-1}\alpha\overline{x}_{i} \in \mathcal{Z}(\mathcal{A}_{n-i})$ and $u_{x_{i}} \in \mathcal{A}_{n-i}.$ This proves that $\sigma$ extends $(\tilde{\sigma}_{H_{i}})_{x_{i}}$ and completes the proof of statement (3.1.1). \qed \para \noindent Theorem \ref{thm4} together with statement (3.1.1) yield the following:
\begin{theorem}{\label{thm5}}
	Let $(H,K)$ be a generalized strong Shoda pair with $H_{i}'s$ and $C_{i}'s$ \linebreak as described above. For each $x_{i} \in C_{i}/H_{i},$ we can always find a unit $u_{x_{i}} \in$ \linebreak $M_{k_{i-1}}(\cdots(M_{k_{0}}(\Q He_{\Q}(\lambda))*^{\sigma_{H_{0}}}_{\tau_{H_{0}}} C_{0}/H_{0}) *^{\sigma_{H_{1}}}_{\tau_{H_{1}}} \cdots)*^{\sigma_{H_{i-1}}}_{\tau_{H_{i-1}}}C_{i-1}/H_{i-1})$ so that
	\begin{enumerate}[(i)]
		\item $z_{x_{i}}=\operatorname{diag}(\overline{x}_{i}u_{x_{i}},\overline{x}_{i}u_{x_{i}},\cdots , \overline{x}_{i}u_{x_{i}})_{k_{i} \cdots k_{n-1}} \in \C_{\mathcal{A}}(\mathcal{B}),$ where $\overline{x}_{i} \in C_{i}$ is a fixed representative of coset $x_{i}.$ \item $z_{x_{i}}^{-1}\mathbb{E}z_{x_{i}}=\mathbb{E}.$ \end{enumerate}
    Furthermore, if we fix such a $z_{x_{i}}$ for all $x_{i} \in C_{i}/H_{i},~0 \leq i \leq n-1,$ then $\{ z_{x_{0}}z_{x_{1}}\cdots z_{x_{n-1}}~|~ x_{i} \in C_{i}/H_{i}, 0 \leq i \leq n-1 \}$ is a basis unit for the crossed product algebra $\operatorname{Cen}_{\mathcal{A}}(\mathcal{B}) \cong (\mathbb{E}/\mathbb{F},\tau).$ Consequently, the Galois group $\G(\mathbb{E}/\mathbb{F})$ and the twisting $\tau$ which appeared in Theorem \ref{thm4} are given by	\begin{enumerate}[(a)]
		\item $\G(\mathbb{E}/\mathbb{F}) = \{\sigma_{x_{0}} \circ \sigma_{x_{1}} \circ \cdots \circ \sigma_{x_{n-1}} ~|~ x_{i} \in C_{i}/H_{i} , ~0 \leq i \leq n-1 \},$ where $\sigma_{x_{i}}$ is the automorphism on $\mathbb{E}$ given by conjugation by $z_{x_{i}}.$
		\item $\tau : \G(\mathbb{E}/\mathbb{F}) \times \G(\mathbb{E}/\mathbb{F}) \rightarrow \mathcal{U}(\mathbb{E})$ sends $(\sigma_{x_{0}} \circ \cdots \circ \sigma_{x_{n-1}},\sigma_{y_{0}} \circ \cdots \circ \sigma_{y_{n-1}})$ to $(z_{x_{0}y_{0}} \cdots z_{x_{n-1}y_{n-1}})^{-1} z_{x_{0}} \cdots z_{x_{n-1}} z_{y_{0}} \cdots z_{y_{n-1}}.$
	\end{enumerate}
\end{theorem}
\subsection{The Particular case, when $\G(\mathbb{E}/\mathbb{F})$ is cyclic}
We now restrict ourselves to the case when $\G(\mathbb{E}/\mathbb{F})$ is cyclic. This happens for example when $[H:K]=1,2,4,p^n,2p^n,p$ an odd prime. 
In this case, $C_{i}/H_{i}$ is also cyclic for any i, $0 \leq i \leq n-1,$ as it is isomorphic to a subquotient of $\G(\mathbb{E}/\mathbb{F})$. If $C_{i}/H_{i}$ is generated by $x_{i}, ~0 \leq i \leq n-1,$ then one can see that
\begin{equation}{\label{eqa}}
\G(\mathbb{E}/\mathbb{F}) = \{ \sigma_{x_{0}}^{m_{0}}\circ \cdots \circ \sigma_{x_{n-1}}^{m_{n-1}}~|~ 0 \leq m_{i} \leq |C_{i}/H_{i}| \}.
\end{equation} Here, $\sigma_{x_{i}},~x_{i} \in C_{i}/H_{i}$ are as defined in section 3.1. Now if $\sigma=\sigma_{x_{0}}^{m_{0}} \circ \cdots \circ \sigma_{x_{n-1}}^{m_{n-1}}$ generates $\G(\mathbb{E}/\mathbb{F})$ for some $m_{i},~0 \leq m_{i} \leq |C_{i}/H_{i}|,$ then $\Q Ge_{\Q}(\lambda^G)$ is seen to be isomorphic to $k \times k$ matrices over the cyclic algebra $(\mathbb{E}/\mathbb{F},\sigma,\textbf{a}),$ where $\textbf{a}= (z_{x_{0}}^{m_{0}} \cdots z_{x_{n-1}}^{m_{n-1}})^{\prod_{i=0}^{n-1}|C_{i}/H_{i}|}.$ In practice, the computation of $\prod_{i=0}^{n-1}|C_{i}/H_{i}|$ power of $z=(z_{x_{0}}^{m_{0}} \cdots z_{x_{n-1}}^{m_{n-1}})$ for a given group $G$ is quite technical but can be computed precisely from the given relations of $G$. One has to succesively compute $z^{|C_{n-1}/H_{n-1}|},$ $z^{|C_{n-1}/H_{n-1}||C_{n-2}/H_{n-2}|},\cdots,$ and finally $z^{\prod_{i=0}^{n-1} |C_{i}/H_{i}|},$ which we shall explain now.\para \noindent
\textbf{(1) \underline{Computation of $z^{|C_{n-1}/H_{n-1}|}$}}
\para \noindent
Observe that the cyclic group generated by $\sigma^{|C_{n-1}/H_{n-1}|}$ is a subgroup of $\G(\mathbb{E}/\mathbb{F})$ of order $\prod_{i=0}^{n-2}|C_{i}/H_{i}|.$ As $\G(\mathbb{E}/\mathbb{F})$ is cyclic, it has a unique subgroup of order $\prod_{i=0}^{n-2}|C_{i}/H_{i}|$ namely $\G(\mathbb{E}/\mathbb{E}_{n-1})$, where $\mathbb{E}_{n-1}$ is defined in equation (\ref{eq}). Hence, $\langle \sigma^{|C_{n-1}/H_{n-1}|} \rangle=\G(\mathbb{E}/\mathbb{E}_{n-1}).$ This means that $z^{|C_{n-1}/H_{n-1}|} \in \C_{\mathcal{A}}(\mathbb{E}_{n-1})$, which is equal to $M_{k_{n-1}}(\mathcal{A}_{1})$. As $z^{|C_{n-1}/H_{n-1}|}$ also belongs to the $\C_{\mathcal{A}}(\mathcal{B})$, it thus follows that \begin{equation}{\label{eqq}}
z^{|C_{n-1}/H_{n-1}|}= \operatorname{diag}(\alpha,\alpha,\cdots,\alpha)_{k_{n-1}k_{n-2}},
\end{equation} for some $\alpha \in \mathcal{A}_{2}*^{\sigma_{H_{n-2}}}_{\tau_{H_{n-2}}} C_{n-2}/H_{n-2}$. Furthermore, from the isomorphism of $\Q Ge_{\Q}(\lambda^G)$ with $M_{k_{n-1}}(M_{k_{n-2}}(\mathcal{A}_{2}*_{\tau_{H_{n-2}}}^{\sigma_{H_{n-2}}} C_{n-2}/H_{n-2})*_{\tau_{H_{n-1}}}^{\sigma_{H_{n-1}}} C_{n-1}/H_{n-1})$ with its precise map given by the repeated application of (\cite{BK3}, Proposition 2(i), also see equation (\ref{eqn})), one can see that $\alpha=e_{\Q}(\lambda^{H_{n-2}}) z^{|C_{n-1}/H_{n-1}|} e_{\Q}(\lambda^{H_{n-2}})$. Hence, in order to compute $z^{|C_{n-1}/H_{n-1}|}$, it is enough to compute $e_{\Q}(\lambda^{H_{n-2}}) z^{|C_{n-1}/H_{n-1}|} e_{\Q}(\lambda^{H_{n-2}})$, which we shall do now. Remember that $z$ is in the $\C_{\mathcal{A}}(\mathcal{B})$ and hence it commutes with $e_{\Q}(\lambda^{H_{i}})$ for all $i.$ Therefore $\alpha=e_{\Q}(\lambda^{H_{n-2}}) z^{|C_{n-1}/H_{n-1}|} e_{\Q}(\lambda^{H_{n-2}})= z^{|C_{n-1}/H_{n-1}|} e_{\Q}(\lambda^{H_{n-2}}).$ Let us now compute $z^{|C_{n-1}/H_{n-1}|} e_{\Q}(\lambda^{H_{n-2}})$, i.e., $\alpha$.
\para \noindent Denote $x_{n-1}^{m_{n-1}}$ by $a_{0}$ and $x_{n-1}^{-(m_{n-1})}ze_{\Q}(\lambda^{H_{n-1}})$ by $b_{0}$ and write $ze_{\Q}(\lambda^{H_{n-1}})=a_{0}b_{0}.$ As $u_{x_{i}}$ belongs to $\mathcal{A}_{n-i},$ we see that $b_{0} \in \mathcal{A}_{1}.$\para \noindent 
Find the least positive integer $d_{0}$ such that 
\begin{flalign}{\label{eq4}}
a_{0}^{-d_{0}} e_{\Q}(\lambda^{H_{n-2}}) a_{0}^{d_{0}} = s^{-1} e_{\Q}(\lambda^{H_{n-2}}) s
\end{flalign}
for some $s \in T_{n-1},$ where $T_{n-1}$ is a transversal of $C_{n-2}$ in $H_{n-1}.$ Such a \linebreak $d_{0}$ exists and it divides $|C_{n-1}/H_{n-1}|.$ The existence of $d_{0}$ follows as \linebreak $a_{0}^{-|C_{n-1}/H_{n-1}|} e_{\Q}(\lambda^{H_{n-2}}) a_{0}^{|C_{n-1}/H_{n-1}|}$ is a $H_{n-1}$-conjugate of $e_{\Q}(\lambda^{H_{n-2}})$ and \linebreak $s^{-1}e_{\Q}(\lambda^{H_{n-2}})s,~ s \in T_{n-1}$ are precisely all the $H_{n-1}$-conjugates of $e_{\Q}(\lambda^{H_{n-2}}).$\para \noindent Let us now show that $d_{0}$ divides $|C_{n-1}/H_{n-1}|.$ Write $|C_{n-1}/H_{n-1}|= d_{0}q+r,$ where $q,r$ are integers and $0 \leq r < d.$ As both $a_{0}^{-|C_{n-1}/H_{n-1}|} e_{\Q}(\lambda^{H_{n-2}}) a_{0}^{|C_{n-1}/H_{n-1}|}$ and $a_{0}^{-d_{0}} e_{\Q}(\lambda^{H_{n-2}}) a_{0}^{d_{0}}$ are $H_{n-1}$-conjugates of  $e_{\Q}(\lambda^{H_{n-2}}),$ it turns out that $a_{0}^{-r} e_{\Q}(\lambda^{H_{n-2}})a_{0}^{r}$ is also a $H_{n-1}$-conjugate of $e_{\Q}(\lambda^{H_{n-2}}),$ which contradicts the minimality of $d_{0}$ unless $r=0.$
\para \noindent In order to find $z^{|C_{n-1}/H_{n-1}|} e_{\Q}(\lambda^{H_{n-2}}),$ we will first compute \linebreak $z^{d_{0}} e_{\Q}(\lambda^{H_{n-2}})$ and then raise it to $\frac{|C_{n-1}/H_{n-1}|}{d_{0}}$ power.
\begin{flalign}{\label{eq5}} \nonumber
z^{d_{0}}e_{\Q}(\lambda^{H_{n-1}}) &= a_{0}^{d_{0}} (a_{0}^{-(d_{0}-1)}b_{0}a_{0}^{(d_{0}-1)} \cdots a_{0}^{-1}b_{0} a_{0} b_{0})\\ 
&=a_{0}^{d_{0}}X,
\end{flalign}
where $X = a_{0}^{-(d_{0}-1)}b_{0}a_{0}^{(d_{0}-1)} \cdots a_{0}^{-1}b_{0} a_{0} b_{0}.$ Observe that $z^{d_{0}}e_{\Q}(\lambda^{H_{n-1}})$ commutes with $e_{\Q}(\lambda^{H_{n-2}})$ as $z$ commutes with $e_{\Q}(\lambda^{H_{i}})$ for all $i$. Therefore,
\begin{equation*}
a_{0}^{d_{0}}Xe_{\Q}(\lambda^{H_{n-2}})=e_{\Q}(\lambda^{H_{n-2}})a_{0}^{d_{0}}X.
\end{equation*} This implies
\begin{equation}{\label{eq6}}
\begin{array}{ll}
Xe_{\Q}(\lambda^{H_{n-2}})X^{-1} &= a_{0}^{-d_{0}}e_{\Q}(\lambda^{H_{n-2}})a_{0}^{d_{0}}\\
&= s^{-1}e_{\Q}(\lambda^{H_{n-2}}) s,~\mbox{by equation}~(\ref{eq4}).
\end{array}
\end{equation}
\noindent The above equation (\ref{eq6}) implies that $sX$ commutes with $e_{\Q}(\lambda^{H_{n-2}}).$ Therefore,$$sX= \left(
\begin{array}{c | c} 
Y & \begin{array}{c c c} 
0 & \cdots & 0 
\end{array} \\ 
\hline 
\begin{array}{c c c} 
0\\
\vdots\\
0 
\end{array} & Z 
\end{array} 
\right) $$ for some $Y \in \mathcal{A}_{2}*^{\sigma_{H_{n-2}}}_{\tau_{H_{n-2}}} C_{n-2}/H_{n-2}$ and $Z$ is a matrix of size $k_{n-2}-1 \times k_{n-2}-1$ over $\mathcal{A}_{2}*^{\sigma_{H_{n-2}}}_{\tau_{H_{n-2}}} C_{n-2}/H_{n-2}.$ The above expression of $sX$ gives that $e_{\Q}(\lambda^{H_{n-2}})sXe_{\Q}(\lambda^{H_{n-2}})=Y,$ which also implies $sXe_{\Q}(\lambda^{H_{n-2}})=Y$ as $sX$ commutes with $e_{\Q}(\lambda^{H_{n-2}})$. \para \noindent
Hence
\begin{flalign*}
z^{d_{0}}e_{\Q}(\lambda^{H_{n-2}}) &=
z^{d_{0}}e_{\Q}(\lambda^{H_{n-1}})e_{\Q}(\lambda^{H_{n-2}})\\
&=a_{0}^{d_{0}}Xe_{\Q}(\lambda^{H_{n-2}}),~~~~~~~~~~~~~\mbox{by equation (\ref{eq5})}\\
&= a_{0}^{d_{0}} s^{-1} sXe_{\Q}(\lambda^{H_{n-2}}),\\
&= a_{0}^{d_{0}} s^{-1}sXe_{\Q}(\lambda^{H_{n-2}}),~~~~~~~\mbox{by equation (\ref{eq6})}\\
&= a_{0}^{d_{0}} s^{-1} Y.
\end{flalign*}
Consequently,
\begin{flalign*}
z^{|C_{n-1}/H_{n-1}|} e_{\Q}(\lambda^{H_{n-2}})&=(z^{d_{0}}e_{\Q}(\lambda^{H_{n-2}}))^{\frac{|C_{n-1}/H_{n-1}|}{d_{0}}}\\
&= (a_{0}^{d_{0}}s^{-1}Y)^{\frac{|C_{n-1}/H_{n-1}|}{d_{0}}}\\
&= (a_{0}^{d_{0}}s^{-1})^{\frac{|C_{n-1}/H_{n-1}|}{d_{0}}}\prod_{k=0}^{{\frac{|C_{n-1}/H_{n-1}|}{d_{0}}}-1} (a_{0}^{d_{0}}s^{-1})^{-k} Y(a_{0}^{d_{0}}s^{-1})^{k}\\
&= a_{1}b_{1}, \mbox{say},
\end{flalign*}
where 
\begin{flalign*}
a_{1}& =(a_{0}^{d_{0}}s^{-1})^{\frac{|C_{n-1}/H_{n-1}|}{d_{0}}},\\
b_{1}&=\displaystyle{\prod_{k=0}^{{\frac{|C_{n-1}/H_{n-1}|}{d_{0}}}-1}} (a_{0}^{d_{0}}s^{-1})^{-k} Y(a_{0}^{d_{0}}s^{-1})^{k}.
\end{flalign*} Finally, substituting the above expression of $z^{|C_{n-1}/H_{n-1}|}e_{\Q}(\lambda^{H_{n-2}})$ in equation (\ref{eqq}), we obtain the following expression of $z^{|C_{n-1}/H_{n-1}|}$;
$$z^{|C_{n-1}/H_{n-1}|}= \operatorname{diag}(a_{1}b_{1},a_{1}b_{1},\cdots, a_{1}b_{1})_{k_{n-1}k_{n-2}}.$$
\para \noindent
\textbf{(2) \underline{Computation of $z^{|C_{n-1}/H_{n-1}||C_{n-2}/H_{n-2}|}$}}
\para\noindent We are basically replacing $a_{0}b_{0}$ with $a_{1}b_{1}$ and following the previous computation. First of all, we will show that $a_{1} \in C_{n-2}.$ 
\begin{flalign*}
H_{n-1}a_{1} &= H_{n-1}(a_{0}^{d_{0}} s^{-1})^{\frac{|C_{n-1}/H_{n-1}|}{d_{0}}}\\
&= H_{n-1} a_{0}^{|C_{n-1}/H_{n-1}|} \hspace{-0.5cm} \prod_{k=0}^{\frac{|C_{n-1}/H_{n-1}|}{d_{0}}-1}\hspace{-0.5cm}(a_{0}^{-d_{0}k}s^{-1}a_{0}^{d_{0}k}).
\end{flalign*}
Observe that \hspace{-0.5cm} $\displaystyle{\prod_{k=0}^{\frac{|C_{n-1}/H_{n-1}|}{d_{0}}-1}}\hspace{-0.5cm}(a_{0}^{-d_{0}k}s^{-1}a_{0}^{d_{0}k})$ belongs to $H_{n-1}$ as $H_{n-1} \unlhd C_{n-1},~s \in H_{n-1}$ and $a_{0} \in C_{n-1}.$ Also $a_{0}^{|C_{n-1}/H_{n-1}|}$ is an element of $H_{n-1}.$ So, the right hand side of above equation is $H_{n-1}$ only. Consequently, we obtain that $a_{1} \in H_{n-1}.$ The commutativity of $a_{1}$ with $e_{\Q}(\lambda^{H_{n-2}})$, by equation (\ref{eq4}), implies that $a_{1} \in C_{n-2}.$ \para \noindent
Analogus to the previous computation, we find the least positive integer $d_{1}$ such that $a_{1}^{-d_{1}} e_{\Q}(\lambda^{H_{n-3}})a_{1}^{d_{1}}= t^{-1}e_{\Q}(\lambda^{H_{n-3}})t$ for some $t \in T_{n-2}.$ Similar arguments reveal that such a $d_{1}$ exist, it divides $|C_{n-2}/H_{n-2}|$ and we finally obtain that $$z^{|C_{n-1}/H_{n-1}||C_{n-2}/H_{n-2}|}= \operatorname{diag}(a_{2}b_{2},a_{2}b_{2},\cdots,a_{2}b_{2})_{k_{n-1}k_{n-2}k_{n-3}},$$
where
\begin{flalign*}
a_{2}& =((a_{1})^{d_{1}}t^{-1})^{\frac{|C_{n-2}/H_{n-2}|}{d_{1}}},\\
b_{2}&= \displaystyle{\prod_{k=0}^{\frac{|C_{n-2}/H_{n-2}|}{d_{1}}-1}} ((a_{1})^{d_{1}}t^{-1})^{-k} t (\prod_{j=0}^{d_{1}-1} a_{1}^{-j}b_{1}a_{1}^j)((a_{1})^{d_{1}}t^{-1})^{k}e_{\Q}(\lambda^{H_{n-3}}).
\end{flalign*} 
\noindent
Repeating the above mentioned arguments, we can succesively compute $z^{|C_{n-1}/H_{n-1}|},$ $z^{|C_{n-1}/H_{n-1}||C_{n-2}/H_{n-2}|}, \cdots,$ and finally, $z^{\prod_{i=0}^{n-1}|C_{i}/H_{i}|}$.
\noindent
In the next two sections, we will apply the theory which has been established in this section to specific examples.
\section{$G = ((C_{p} \times C_{p}) \rtimes C_{p}) \rtimes C_{2^n}$}
Let $p$ be an odd prime, $p \equiv 1 \operatorname{mod} 4$ and let $$P= \langle x,y,z~|~ x^{p}=y^{p}=z^{p}=e, xy=yx, yz=zyx,xz=zx\rangle $$ be the extraspecial $p$-group of order $p^3.$ Write $p-1 = 2^{n-1}t, ~t$ odd and $n \geq 3$. Let $r$ be a primitive root $\operatorname{modulo} p$, i.e., the order of $r$ modulo $p$ is $p-1$. Let $k,q$ be integers with $r^{t} \equiv -kq \operatorname{mod} p.$ Consider the automorphism $\theta: P \rightarrow P$ given by $x \mapsto x^{r^{t}}, y \mapsto z^{k}, z \mapsto y^{q}.$ It can be checked that the order of $\theta$ is $2^n.$ Consider $P \rtimes C_{2^n},$ the semidirect product of $P$ with the cyclic group of order $2^n$. It has generators $x,y,z,w$ with the following relations:
\begin{quote}
	$x^{p}=y^{p}=z^{p}=w^{2^n}=e, xy=yx, yz=zyx, xz=zx, w^{-1}xw=\theta(x),\linebreak w^{-1}yw=\theta(y),
	w^{-1}zw=\theta(z).$
\end{quote}
Let us denote the group $P \rtimes C_{2^n},$ by $G_{p^32^n}.$ In this section, we will show that $ G_{p^32^n}$ is a generalized strongly monomial group and will illustrate the theory developed in previous section to write the simple components along with their Schur indices and the precise Wedderburn decomposition of $\Q  G_{p^32^n}.$

\subsection{Irreducible characters}
Denote by $Q$ the subgroup of $G_{p^32^n}$ generated by $x,y,w^{2^{n-1}}.$

\begin{theorem}{\label{thm6}}
	If $G=G_{p^32^n},$ the following statements hold:
	\begin{enumerate}[(i)]
		\item there are $2^n$ linear characters of $G$ and these are given by $\lambda_{i}: G     \rightarrow \mathbb{C}$ such that $x \mapsto 1, y \mapsto 1, z \mapsto 1, w \mapsto \zeta^{i},$ where $0 \leq i \leq 2^n-1$ and $\zeta$ is a primitive $2^n$-th root of unity in the field $\mathbb{C}$ of complex numbers.
		\item every non-linear irreducible character of $G$ is either \\
		(a) $\phi^{G}$ for some non-principal linear character $\phi$ of $P$ or\\
		(b) $\psi^{G}$ for some non-principal linear character $\psi$ of $Q,$ such that $x \notin \operatorname{Ker}\psi.$	
	\end{enumerate}
\end{theorem}

\begin{lemma}{\label{lem2}}
	If $G=G_{p^32^n}$ and $\phi$ is a non-principal linear character of $P,$ then $\phi^{G}$ is irreducible. This way one obtains $\frac{p^2-1}{2^n}$ distinct irreducible characters of $G$ which are induced from the linear characters of $P.$
\end{lemma}
\noindent {\bf Proof.} Suppose $\phi$ is a non-principal linear character of $P.$ As the commutator subgroup of $P$ is $\langle x \rangle$, the kernel of $\phi$ is a normal subgroup of $P$ with cyclic quotient and it contains $\langle x \rangle$. This gives that the $\operatorname{ker}(\phi)$ is one of $\langle x,z\rangle , \langle x,z^{i}y\rangle , 0 \leq i \leq p-1.$ In view of (\cite{Issac}, Exercise 6.1), to show that $\phi^{G}$ is irreducible, it is enough to show that $I_{G}(\phi) = P,$ where $I_{G}(\phi)= \{g \in G~|~\phi^g=\phi, i.e., \phi(g^{-1} \alpha g) = \phi(\alpha)~\forall~\alpha \in P \}$ is the inertia of $\phi$ in $G.$ Let us assume that $w^{j} \in I_{G}(\phi)$ for some $0 \leq j \leq 2^n-1.$ It means $\phi(w^{-j} \alpha w^j) = \phi(\alpha)~ \forall ~\alpha \in P,$ i.e., $w^{-j}\alpha w^j\alpha^{-1}$ belongs to $\operatorname{Ker}\phi$ for all $\alpha$ in $P.$ Suppose $\operatorname{Ker}\phi = \langle x,z \rangle.$ Then $y \notin \operatorname{Ker}\phi.$ This is a contradiction, as from the group relations, we have
$$ w^{-j}\alpha w^j \alpha^{-1} = \left\{
                                    \begin{array}{ll}
                                    y^{(kq)^{j/2}-1}, & \mbox{where} ~\alpha=y,~ j~ \mbox{even};\\
                                    y^{k^{\frac{j-1}{2}} q^{\frac{j+1}{2}}}z^{-1}, & \mbox{where}~ \alpha = z,~j~\mbox{odd}.
                                    \end{array}
\right.
$$ Thus we obtain that $w^j$ does not belong to $\operatorname{I}_{G}(\phi)$ for any $1 \leq j \leq 2^n-1$. Hence, $\operatorname{I}_{G}(\phi)=P$ giving that $\phi^G$ is irreducible. \para \noindent
Applying the similar arguments, when $\operatorname{Ker}\phi= \langle x,y \rangle$ or $\langle x, z^iy \rangle,~1 \leq i \leq p-1,$ it can be checked that $w^j \notin \operatorname{I}_{G}(\phi)$ for any $j$. Therefore in these cases also, $\operatorname{I}_{G}(\phi)=P$, i.e., $\phi^G$ is irreducible. \para \noindent The last statement of the lemma follows because there are total $p^2-1$ non-principal linear characters of $P,$ and if $\phi, ~\phi'$ are two linear characters of $P,$ then $\phi^G = \phi'^G$ if and only if $\phi' = \phi^{w^{j}}$ for some $j,~0 \leq j \leq 2^n-1.$\qed	

\begin{lemma}{\label{lem3}}
	If $G = G_{p^32^n}$ and $\psi$ is a non-principal linear character of $Q$ such that $x \notin \operatorname{Ker}\psi,$ then $\psi^{G}$ is irreducible. This way one obtains $\frac{p-1}{2^{n-2}}$ distinct irreducible characters of $G$ which are induced from linear characters of $Q.$
\end{lemma}\noindent
\textbf{Proof.} Suppose $\psi$ is a non-principal linear character such that $x \notin \operatorname{Ker}\psi.$ Let $H_{1}$ be a subgroup of $G$ generated by $x,y,z,w^{2^{n-1}}.$ To show that $\psi^{G}$ is irreducible, we will show the following:
\begin{enumerate}[(i)]
	\item $\psi^{H_{1}}$ is irreducible;
	\item $(\psi^{H_{1}})^G$ is irreducible.
\end{enumerate}
For (i), in view of Shoda's Theorem, (see \cite{shoda}), it is enough to prove that for every $g \in H_{1} \backslash Q,$ there exists $h \in Q \cap Q^g$ such that $\psi(ghg^{-1}) \neq \psi(h).$ Suppose $g \in H_{1} \backslash Q.$ Write $g = qz^i,$ where $q \in Q$ and $0 < i \leq p-1.$ Observe that $y \in Q \cap Q^g$ and $\psi(gyg^{-1})=\psi(qz^iyz^{-i}q^{-1}) = \psi(qyx^{-i}q^{-1}),$ which is never equal to $\psi(y)$, as $qyx^{-i}q^{-1}$ is either $yx^{-i}$ or $y^{-1}x^{-i}$ and $x^{i} \notin \operatorname{Ker}\psi$ for $i > 0.$ Hence, (i) follows.\para \noindent
To prove (ii), as $H_{1} \unlhd~ G,$ it is enough to prove that $I_{G}(\psi^{H_{1}})= H_{1}$ (see \cite{Issac}, Exercise 6.1). For this, it only suffices to prove that $w^j,$ for any $0 < j \leq 2^{n-1}-1,$ does not belong to $I_{G}(\psi^{H_{1}}).$ Suppose on the contrary, $w^j \in I_{G}(\psi^{H_{1}})$ for some $0 < j \leq 2^{n-1}-1$. This gives $(\psi^{H_{1}})^{w^j}(x)= \psi^{H_{1}}(x).$ As $(\psi^{H_{1}})^{w^{j}}(x)= p\psi(x^{r^{tj}})$ and $\psi^{H_{1}}(x)= p\psi(x),$ we obtain that $\psi(x^{r^{tj}}) = \psi(x).$ Since $x \notin \operatorname{Ker}\psi,~\psi(x)$ is a root of unity of order $p.$ Hence $r^{tj} \equiv 1 (\operatorname{mod} p),$ i.e.,  $p-1$ divides $tj.$ It is a contradiction, as $p-1= 2^{n-1}t$ and $j < 2^{n-1}.$ Hence (ii) is proved.\para \noindent We are now going to show that by inducing to $G$ the linear characters of $Q$ which do not contain $x$ in its kernel, we obtain $\frac{p-1}{2^{n-2}}$ distinct irreducible characters of $G$. Observe that the only normal subgroups $N$ of $Q$ such that $Q/N$ is cyclic and $x \notin \operatorname{ker}\psi$ are $N=\langle y \rangle$ or $\langle y,w^{2^{n-1}} \rangle$. Hence we obtain that the kernel of $\psi$ is one of $\langle y \rangle$ or $\langle y,w^{2^{n-1}} \rangle$. If $\operatorname{Ker}\psi = \langle y \rangle,$ then there are total $p-1$ such linear characters $\psi_{i}: Q \rightarrow \mathbb{C}$ given by $x \mapsto \zeta_{p}^{i},~y \mapsto 1,~w^{2^{n-1}} \mapsto -1,~1 \leq i \leq p-1.$ As, $\psi_{i}^{H_{1}}(x)= p\zeta_{p}^{i},$ it follows that all $\psi_{i}^{H_{1}},~1 \leq i \leq p-1,$ are distinct irreducible characters of $H_{1}.$ Out of these $p-1$ characters, $\frac{p-1}{2^{n-1}}$ will be distinct when induced to $G$ because $\operatorname{I}_{G}(\psi_{i}^{H_{1}})=H_{1}$. This gives $\frac{p-1}{2^{n-1}}$ irreducible characters of $G,$ which are induced from linear character $\psi$ of $Q$ with $\operatorname{Ker}\psi= \langle y \rangle.$ Similarly, we have $\frac{p-1}{2^{n-1}}$ irreducible characters of $G$ which are induced from linear characters $\psi$ of $Q$ with $\operatorname{Ker}\psi= \langle y, w^{2^{n-1}} \rangle.$ Consequently, we have $\frac{p-1}{2^{n-2}}$ distinct irreducible characters of $G$ which are induced from linear characters of $Q.$ This proves the theorem. \qed \para \noindent
\textbf{Proof of Theorem \ref{thm6}.} The statement (i) follows because the commutator subgroup of $G$ is $P$. Further, the number of linear characters of $G$ is $|G/P|$ which is $2^n.$ From Lemmas \ref{lem2} and \ref{lem3}, we have obtained $\frac{p^2-1}{2^n}$ distinct irreducible characters of degree $2^n$ and $\frac{p-1}{2^{n-2}}$ distinct irreducible characters of degree $p.2^{n-1}.$ As $$2^{n}.1^2+\frac{p^2-1}{2^n}.2^{2n}+\frac{p-1}{2^{n-2}}.p^2.2^{2n-2}=2^n.p^3 = |G|,$$ it turns out that we have obtained all the irreducible characters of $G$. Hence, (ii) is proved. \qed
\subsection{Shoda Pairs}
In Theorem {\ref{thm6}}, we have seen that all characters of $G=G_{p^32^n}$ are monomial and hence it is a monomial group. Let us now see the Shoda pairs of $G$ which arise from these monomial irreducible characters. Clearly, the linear characters yield the following Shoda pairs of $G:$
\begin{equation}{\label{equ7}}
(G,\langle P,w^{2^{i}}\rangle ), 0 \leq i \leq n.
\end{equation}
Let us determine the Shoda pairs arising from the irreducible characters of $G$ induced from $P.$ Observe that the only possible kernels of linear characters of $P$ are $\langle x,z\rangle , \langle x,z^jy\rangle ,~0 \leq j \leq p-1,$ as these are precisely all the subgroups of $P$ with cyclic quotient. So, irreducible characters which are induced from $P$ yield the following Shoda pairs :
\begin{equation}{\label{equ8}}
(P,\langle x,z\rangle ), (P,\langle x,z^jy\rangle ), 0 \leq j \leq p-1.
\end{equation}
Finally, one can check that the following Shoda pairs are obtained from the irreducible characters of $G$ which are induced from $Q:$
\begin{equation}{\label{equ9}}
(Q,\langle y\rangle ), (Q,\langle y,w^{2^{n-1}}\rangle ).
\end{equation}
Let us now identify the inequivalent Shoda pairs which we have obtained in equations $(\ref{equ7}),(\ref{equ8})$ and $(\ref{equ9}).$ \para \noindent Recall from Remark \ref{remark1}(b), we have that the Shoda pairs which arise from irreducible characters of distinct degree cannot be equivalent. Hence, no Shoda pair in equation $(i)$ can be equivalent to that in equation $(j)$, $i,j= \ref{equ7},\ref{equ8},\ref{equ9}$ and $i \neq j.$ \para \noindent Also from Remark \ref{remark1}(a), we have that if $H \unlhd G,$ then the Shoda pair $(H,K_{1})$ is equivalent to $(H,K_{2})$ if and only if $K_{1}$ and $K_{2}$ are conjugate in $G.$ Applying this, we obtain that the Shoda pair given in equation (\ref{equ7}) are all inequivalent as $\langle P,w^{2^{i}}\rangle ,~0 \leq i \leq n,$ are all normal in $G.$\para \noindent Next, the only $G$-conjugate of $\langle x,z\rangle $ is $\langle x,y\rangle $ and that of $\langle x,z^jy\rangle $ is $\langle x,z^{kq'j'}y\rangle,$ where $qq' \equiv 1 \operatorname{mod} p$ and $jj' \equiv 1 \operatorname{mod} p.$ Therefore, the following inequivalent Shoda pairs are obtained from equation $(\ref{equ8}):$
\begin{equation}{\label{eqqq}}
	(P,\langle x,z^jy\rangle ), ~j \in J,
\end{equation} 	where $J=\{0\} \cup J'$ and $J'$ is a subset of $\{1,2,\cdots,p-1\}$ of cardinality $\frac{p-1}{2}$ consisting of one representative from the $G$-conjugates of $\langle x,z^jy \rangle,~1\leq j \leq p-1.$   \para \noindent
Moving to the Shoda pairs in equation $(\ref{equ9}),$ we assert that both of them are inequivalent. Consider the linear characters $\lambda_{1}$ and $\lambda_{2}$ on $Q$ given by $\lambda_{1}(x)=\zeta_{p} ,~\lambda_{1}(y) =1,~ \lambda_{1}(w^{2^{n-1}})=-1$ and  $\lambda_{2}(x)=\zeta_{p} , ~\lambda_{2}(y) =1, ~\lambda_{2}(w^{2^{n-1}})=1$, where $\zeta_{p}$ is primitive $p$-th root of unity. If $(Q,\langle y\rangle )$ and $(Q,\langle y,w^{2^{n-1}}\rangle )$ are equivalent, then $e_{\Q}(\lambda_{2}^{G}) = e_{\Q}(\lambda_{1}^{G}),$ i.e., $\lambda_{1}^{G} = \sigma \circ \lambda_{2}^G,$ for some automorphism $\sigma$ of $\Q(\lambda_{2}^{G}).$ Evaluting at $w^{2^{n-1}},$ we obtain that
$$ \lambda_{1}^{G}(w^{2^{n-1}}) = \sigma(\lambda_{2}^{G}(w^{2^{n-1}})).$$   
However, $\lambda_{1}^{G}(w^{2^{n-1}}) = -2^{n-1}$ and $\lambda_{2}^{G}(w^{2^{n-1}}) = 2^{n-1},$ which gives that $\sigma(2^{n-1})=-2^{n-1},$ i.e., $\sigma(1)=-1,$ which is not true. This proves the assertion. \para \noindent
We now want to see that the Shoda pairs which we have obtained are all generalized strong Shoda pairs. From Remark {\ref{remark2}$(ii)$}, it turns out that the Shoda pairs in equations $(\ref{equ7})$ and $(\ref{equ8})$ are strong Shoda pairs. Let us now show that the two Shoda pairs given in equation $(\ref{equ9})$ are generalized strong Shoda pairs of G. \para \noindent Suppose $(H,K)= (Q,\langle y\rangle ).$ We assert that $$Q=\langle x,y,w^{2^{n-1}}\rangle  \leq H_{1}= \langle x,y,z,w^{2^{n-1}}\rangle  \leq G=\langle x,y,z,w\rangle $$ is a strong inductive chain. As $H_{1} \unlhd G,$ in view of Remark \ref{remark2}$(i)$, we only have to see the following: 
\begin{enumerate}[(i)]
	\item $H \unlhd \operatorname{Cen}_{H_{1}}(e_{\Q}(\lambda_{1}));$
	\item the distinct $H_{1}$-conjugates of $e_{\Q}(\lambda_{1})$ are mutually orthogonal.
\end{enumerate}
By (\cite{Jes}, Lemma 3.5.1 (3),(4)), $\operatorname{Cen}_{H_{1}}(e_{\Q}(\lambda_{1})) = N_{H_{1}}(K).$ We see that $N_{H_{1}}(K)= Q,$ as $z$ does not normalize $K.$ Therefore, $\operatorname{Cen}_{H_{1}}(e_{\Q}(\lambda_{1}))$ is $Q$ and hence (i) holds. Let us now show that $e_{\Q}(\lambda_{1}) e_{\Q}(\lambda_{1})^{z^{i}} = 0, ~\mbox{for all} ~ 1 \leq i \leq p-1.$ \para \noindent
Suppose $\lambda_{1}^{'}$ is the restriction of $\lambda_{1}$ on $\langle x,y\rangle .$ By (\cite{BK3}, Lemma 1),
$$e_{\Q}(\lambda_{1}) e_{\Q}(\lambda_{1}^{'})= e_{\Q}(\lambda_{1})= e_{\Q}(\lambda_{1}') e_{\Q}(\lambda_{1})$$ and $$e_{\Q}(\lambda_{1}')^{z^{i}} e_{\Q}(\lambda_{1})^{z^i} =e_{\Q}(\lambda_{1})^{z^i} =  e_{\Q}(\lambda_{1})^{z^i} e_{\Q}(\lambda_{1}')^{z^i}.$$ Therefore, $e_{\Q}(\lambda_{1}) e_{\Q}(\lambda_{1})^{z^{i}} = e_{\Q}(\lambda_{1}) e_{\Q}(\lambda_{1}')  e_{\Q}(\lambda_{1}')^{z^{i}} e_{\Q}(\lambda_{1})^{z^{i}},$ which is zero as $e_{\Q}(\lambda_{1}')$ and $e_{\Q}(\lambda_{1}')^{z^{i}},$ being distinct primitive central idempotents of $\Q \langle x,y\rangle$, are orthogonal. This proves (ii). Hence we have shown that $(H,K)=(Q,\langle y\rangle )$ is a generalized strong Shoda pair of $G.$ \para \noindent The similar arguments can be applied to the Shoda pair $(H,K)=(Q,\langle y,w^{2^{n-1}}\rangle )$ and it turns out that it is also a generalized strong Shoda pair of $G_{p^32^n}$ with the following as strong inductive chain: $$ Q=\langle x,y,w^{2^{n-1}}\rangle  \leq H_{1}= \langle x,y,z,w^{2^{n-1}}\rangle  \leq G=\langle x,y,z,w\rangle.$$ The above discussion yield the following:

\begin{theorem}
	The group $G= G_{p^32^n}$ is a generalized strongly monomial group and it has a complete and irreduntant set of Shoda pairs consisting of the following:
	\begin{enumerate}[(i)]
		\item $(G,\langle P,w^{2^{i}}\rangle ), 0 \leq i \leq n;$
		\item $(P,\langle x,z^jy\rangle ), j \in J;$
		\item $(Q,\langle y\rangle ), (Q,\langle y,w^{2^{n-1}}\rangle ),$
	\end{enumerate}
	where $J$ has cardinality $\frac{p+1}{2}$ and is as defined in equation (\ref{eqqq}).
\end{theorem}

\subsection{Simple Components}  
The Shoda pairs $(G, \langle P,w^{2^{i}}\rangle),~0 \leq i \leq n$ and $(P,\langle x,z^jy\rangle ),~j \in J$, are strong Shoda pairs of $G=G_{p^32^n}$. If $(H,K)= (G,\langle P,w^{2^{i}}\rangle ),$ then by (\cite{Jes}, Lemma 3.3.2), the simple component of $\Q G$ realized by $(H,K)$ is isomorphic to  $\Q(\zeta_{2^{n-i}})$,  where $\zeta_{2^{n-i}}$ is a primitive $2^{n-i}$-th root of unity. Next suppose $(H,K)= (P, \langle x,z^jy \rangle)$ for some $j \in J.$ As the normalizer of $\langle x, z^jy \rangle$ in $G$ is $\langle P, {w^2} \rangle,$ it follows from (\cite{Jes}, Theorem 3.5.5, also see Theorem \ref{thm2}) that the corresponding simple component  is isomorphic to $M_{2}(\Q(\zeta_{p}) *^{\sigma}_{\tau} \langle P,w^{2} \rangle / P),$ where $\sigma_{w^2}$ sends $\zeta_{p}$ to $\zeta_{p}^{-r^t}$ and the twisting $\tau$ is trivial. The triviality of the twisting $\tau$ and the faithful action $\sigma$ implies, by (\cite{Rei}, Corollary 29.8), that the algebra $ \Q(\zeta_{p})*^{\sigma}_{\tau} \langle P,w^{2}\rangle /P $ splits, and consequently the corresponding simple component is isomorphic to $M_{2^n}(F),$ where $F$ is the fixed subfield of $\Q(\zeta_{p})$ which is kept fixed by $\langle \sigma_{w^2} \rangle.$ As  $(\zeta_{p} + \zeta_{p}^{r^t}+\cdots + \zeta_{p}^{r^{(2^{n-1}-1)t}})$ is fixed by $\sigma_{w^2}$ and $[\Q(\zeta_{p}):\Q(\zeta_{p} + \zeta_{p}^{r^t}+\cdots + \zeta_{p}^{r^{(2^{n-1}-1)t}})]= 2^{n-1},$ one obtains that $F=\Q(\zeta_{p} + \zeta_{p}^{r^t}+\cdots + \zeta_{p}^{r^{(2^{n-1}-1)t}}).$ Hence, we have seen that
$$\Q Ge_{\Q}(\lambda^G) \cong \left\{
\begin{array}{ll}
\Q(\zeta_{2^{n-i}})~ if~ (H,K)=(G,\langle P,w^{2^{i}}\rangle ),~0 \leq i \leq n,\\
M_{2^{n}}(F)~ if~ (H,K)=(P,\langle x,z^jy\rangle )~,j \in J.
\end{array}
\right.
$$
The following result provides the structure of simple component realized by the remaining two Shoda pairs. 
\begin{lemma}
	If $G=G_{p^32^n}$, $(H,K)= (Q,\langle y\rangle )$ or $(Q,\langle y,w^{2^{n-1}}\rangle )$ and $\lambda$ a linear characters on $H$ with kernel $K,$  then
	$$
	\Q Ge_{\Q}(\lambda^G) \cong \left\{
	\begin{array}{ll}
	M_{p2^{n-1}}(\mathbb{H}(F)), & {\rm if}~ (H,K)=(Q,\langle y\rangle );\\
	M_{p2^{n}}(F), & {\rm if}~ (H,K)=(Q,\langle y,w^{2^{n-1}}\rangle ),
	\end{array}
	\right.
	$$
	where $F= \Q (\zeta_{p} + \zeta_{p}^{r^t}+\cdots + \zeta_{p}^{r^{(2^{n-1}-1)t}})$ and $\mathbb{H}(F)$ is the quaternion algebra over $F.$
\end{lemma}
\noindent
\textbf{Proof.} We first deal with the simple component realized by $(H,K)=(Q,\langle y\rangle ).$ In section 4.2, we have seen that $Q \leq H_1 \leq G,$ where $H_{1}=\langle x,y,z,w^{2^{n-1}}\rangle,$ is a strong inductive chain for the Shoda pair $(Q,\langle y\rangle).$  Consider $\lambda: Q \rightarrow \mathbb{C}$ given by $x \mapsto \zeta_p, y \mapsto 1, w^{2^{n-1}} \mapsto -1.$ \para \noindent {\underline{Step 1: Compute $C_{0}= \operatorname{Cen}_{H_{1}}(e_{\Q}(\lambda))$ and $C_{1}= \operatorname{Cen}_{G}(e_{\Q}(\lambda^{H_{1}})).$} \para \noindent 
In the proof of Lemma 3, we have already seen that $z \notin \operatorname{Cen}_{H_{1}}(e_{\Q}(\lambda))$ and hence $C_{0}=Q.$ Let us now show that $C_{1}=G.$ For an arbitrary element $x^ay^bz^c(w^{2^{n-1}})^{d} \in H_{1},$ where $0 \leq a,b,c \leq p-1$ and $d= 0,1,$ we observe that
$$ \lambda^{H_{1}}(x^ay^bz^c(w^{2^{n-1}})^{d}) = \left\{
\begin{array}{ll}
-\zeta_p^{a+bc/2} ~~~~~ \mbox{if}~~ d \neq 0;\\
p \zeta_p^a ~~~~~~~~~~~~\mbox{if}~~ d=0,c=0,b=0;\\
0~~~~~~~~~~~~~~~\mbox{if}~~ d=0, c \neq 0~ \mbox{or}~b \neq 0
\end{array}
\right.
$$
and
$$ (\lambda^{H_{1}})^{w}(x^ay^bz^c(w^{2^{n-1}})^{d}) = \left\{
\begin{array}{ll}
-\zeta_p^{ar^t+r^tbc/2} ~~~~~ \mbox{if}~~ d \neq 0;\\
p \zeta_p^{ar^t} ~~~~~~~~~~~~~~\mbox{if}~~ d=0,c=0,b=0;\\
0~~~~~~~~~~~~~~~~~~\mbox{if}~~ d=0, c \neq 0~ \mbox{or}~b \neq 0.
\end{array}
\right.
$$
Hence, if $\sigma$ is the automorphism of $\Q(\zeta_{p})$ which sends $\zeta_{p}$ to $\zeta_{p}^{r^t},$ then $(\lambda^{H_{1}})^{w}= \sigma \circ \lambda^{H_{1}}.$ Therefore, $e_{\Q}(\lambda^{H_{1}})^{w}= e_{\Q}(\lambda^{H_{1}}),$ i.e., $w \in \operatorname{Cen}_{G}(e_{\Q}(\lambda^{H_{1}})).$ This gives, $C_{1}= \operatorname{Cen}_{G}(e_{\Q}(\lambda^{H_{1}}))= G.$ \para \noindent 
{\underline{Step 2:  Use Theorem \ref{thm4}.}} \para \noindent As $(H,K)$ is a strong Shoda pair of $H_{1}$ and $C_{0} = Q, k_{0} =p,$  by (\cite{BK3}, Proposition 2(i), also see equation (\ref{eqn})), \begin{equation}{\label{equ10}} \Q H_{1}e_{\Q}(\lambda^{H_{1}}) \stackrel{\theta}{\cong} M_{p}(\Q(\zeta_{2p})), \end{equation} where  $\theta$ sends  any $\alpha \in \Q H_{1}e_{\Q}(\lambda^{H_{1}}) $ to $(\alpha_{ij})_{p \times p}$, where $\alpha_{ij} = e_{\Q}(\lambda) z^{j-1} \alpha z^{-(i-1)} e_{\Q}(\lambda).$\para \noindent
Further, as $C_{1}=G$ and $k_{1}=1,$ by (\cite{BK3}, Proposition 2(i)),
\begin{flalign*}
\Q Ge_{\Q}(\lambda^G) & \cong M_{k_{1}}(\Q H_{1}e_{\Q}(\lambda^{H_{1}})*^{\sigma_{H_{1}}}_{\tau_{H_{1}}} C_{1}/H_{1}), \\&
= M_{p}(\Q(\zeta_{2p}))*^{\sigma_{H_{1}}}_{\tau_{H_{1}}}  C_{1}/H_{1}.
\end{flalign*}
Using Theorem \ref{thm4}, we have $$ M_{p}(\Q(\zeta_{2p}))*^{\sigma_{H_{1}}}_{\tau_{H_{1}}}  C_{1}/H_{1}\cong M_{p}(\mathbb{E}/\mathbb{F},\tau),$$ where $\mathbb{E}=\{\operatorname{diag}(\alpha,\alpha,\cdots,\alpha)_{p}~|~\alpha \in \Q He_{\Q}(\lambda)\}$ and $\mathbb{F}=\mathcal{Z}(M_{p}(\Q(\zeta_{2p}))*^{\sigma_{H_{1}}}_{\tau_{H_{1}}}  C_{1}/H_{1})$, which is subfield of $\mathbb{E}$ of index $|C_{0}/H_{0}||C_{1}/H_{1}|=2^{n-1}$ by step (ii) of section 3. Also the twisting $\tau$ is as explained in the Theorem \ref{thm5} and we will describe $\tau$ more precisely in the following steps.
\para \noindent {\underline{Step 3:  Compute  $\G(\mathbb{E}/\mathbb{F}) $.}} \para \noindent We have from  section 3.1 that  $\G(\mathbb{E}/\mathbb{F}) $ is of type $C_0/H_0$-by-$C_1/H_1$.  As $C_0= H_0$ and  $C_1 = G$, we obtain  that  $\G(\mathbb{E}/\mathbb{F}) $  is isomorphic to $G/H_1$, the cyclic group generated by $H_{1}w$ of order $2^{n-1}$.  Furthermore, from equation (\ref{eqa}), we  see that $$\G(\mathbb{E}/\mathbb{F} )= \langle \sigma_{w}\rangle, $$
where $\sigma_w$ is the automorphism of $\mathbb{E}$ given by conjugation by some $z_w \in \operatorname{Cen}_{\mathcal{A}}(\mathcal{B})$.  Moreover, we have from (3.1.2)  that $z_w = wA_w$ for  some $A_w \in \Q H_1e_{\Q}(\lambda^{H_1})  \cong M_{p}(\Q(\zeta_{2p})).$\para \noindent {\underline{Step 4:  Compute  $A_w$ and hence $\tau$.}} \para \noindent We wish to find $A_w \in M_{p}(\Q(\zeta_{2p}))$ so that the following are satisfied: 
\begin{enumerate}[(i)]
	\item $wA_{w} \in \operatorname{Cen}_{\mathcal{A}}(\mathcal{B}), $ where  $\mathcal{B}=M_{p}(\mathbb{F})$; 
	\item $z_{w}^{-1}\mathbb{E}z_{w}= \mathbb{E}$ and
	\item the automorphism  $\alpha \mapsto z_{w}^{-1} \alpha z_{w}$ of $\mathbb{E}$ has order $2^{n-1}.$
\end{enumerate} 
Suppose $A_{w}=(a_{ij})_{p \times p}$ where $a_{ij} \in \Q He_{\Q}(\lambda).$ Then $wA_{w} \in \operatorname{Cen}_{\mathcal{A}}(\mathcal{B})$ if and only if 
\begin{equation}{\label{equ12}}
A_{w}^{-1}w^{-1}BwA_{w}=B~ \mbox{for all}~ B \in \mathcal{B}.
\end{equation}
Let us denote $w^{-1}Bw$ by $C$ and write $B=(b_{ij})_{p \times p},~C=(c_{ij})_{p \times p},$ where $b_{ij}$ , $C_{ij} \in \mathbb{F}.$ As $c_{ij}$ is the $(i,j)th$ entry of $\theta(w^{-1}\theta^{-1}(B)w)),$ it turns out that $$c_{ij}=\sum \sum x^{(j-1)(l-1)-(i-1)(m-1)q} b_{ml}  e_{\Q}(\lambda)/p.$$ Consequently, to satisfy equation (\ref{equ12}), we need $A_w$  to be such that 
\begin{equation*}
A_{w}^{-1}(c_{ij}) A_{w} = (b_{ij}),
\end{equation*} where $b_{ij} \in \mathbb{F}$ and  $c_{ij}=\sum \sum x^{(j-1)(l-1)-(i-1)(m-1)q} b_{ml}  e_{\Q}(\lambda)/p, $ $1 \leq i, j \leq p$.  
By taking $a_{ij}= x^{-(i-1)(j-1)q}$ and $A_{w}=(a_{ij})_{p \times p},$  we see that the above equation is satisfied and consequently (i) holds. Also, we have  checked that (ii) and (iii) hold for this $A_{w}.$  \para \noindent
With this choice of $A_{w}$ and correspondingly $z_w = wA_w$, we obtain that the twisting $\tau : \G(\mathbb{E}/\mathbb{F}) \times \G(\mathbb{E}/\mathbb{F}) \rightarrow \mathcal{U}(\mathbb{E})$ is given by:
$$\tau(\sigma^i,\sigma^j) = \left\{
\begin{array}{ll}
1, & {\rm if}~i+j  < 2^{n-1};\\
z_{w}^{2^{n-1}}, & {\rm if}~i+j \geq  2^{n-1}.
\end{array}
\right.
$$ Hence  $\Q Ge_{\Q}(\lambda^G)$ becomes isomorphic to the cyclic algebra  $(\mathbb{E}/\mathbb{F},\sigma_w, z_w^{2^{n-1}}).$
\para \noindent {\underline{Step 5:  Compute  $z_w^{2^{n-1}}.$}} \para \noindent In order to determine the Schur index of the cyclic algebra  $(\mathbb{E}/\mathbb{F},\sigma_w, z_w^{2^{n-1}})$, we need to compute $z_w^{2^{n-1}}$, where  $z_w=wA_w$, $A_{w}=(a_{ij})$ and $a_{ij}=x^{-(i-1)(j-1)q}e_{\Q}(\lambda).$  We will follow the procedure mentioned  in section 3.2 to determine the same. Observe that $d_0 =2$, where $d_0$ is as defined in equation (\ref{eq4}). Now $$ z_{w}^2=wA_{w}wA_{w}=w^2X,$$  where $X=w^{-1}A_{w}wA_{w}.$ As $z_{w}^2$ and $w^2$ both commute with $e_{\Q}(\lambda),$ it turns out that $X$ commutes with $e_{\Q}(\lambda)$ and hence 
$$X= \left(
\begin{array}{c | c} 
pe_{\Q}(\lambda) & \begin{array}{c c c} 
0 & \cdots & 0 
\end{array} \\ 
\hline 
\begin{array}{c c c} 
0\\
\vdots\\
0 
\end{array} & Z 
\end{array} 
\right) $$ where $Z$ is a $p-1 \times p-1$ matrix over $\Q He_{\Q}(\lambda).$ Therefore,	\begin{flalign*}
z_{w}^2 e_{\Q}(\lambda) &=  w^2 X e_{\Q}(\lambda)\\
&= w^2pe_{\Q}(\lambda).
\end{flalign*}
Consequently,
\begin{flalign*}
z_{w}^{2^{n-1}} e_{\Q}(\lambda) &= (z_{w}^2)^{2^{n-2}} e_{\Q}(\lambda)\\
&= (z_{w}^2 e_{\Q}(\lambda))^{2^{n-2}} \\
&= w^{2^{n-1}}p^{2^{n-2}}e_{\Q}(\lambda),
\end{flalign*} Thus, from equation (\ref{eqq}), we have that 	\begin{equation}{\label{equ13}} z_{w}^{2^{n-1}} = \operatorname{diag}(w^{2^{n-1}}p^{2^{n-2}}e_{\Q}(\lambda),w^{2^{n-1}}p^{2^{n-2}}e_{\Q}(\lambda),\cdots,w^{2^{n-1}}p^{2^{n-2}}e_{\Q}(\lambda))_{p}\end{equation}
\para \noindent	{\underline{Step 6:  Finally, identify  $\Q Ge_{\Q}(\lambda^G).$}}\para \noindent We have $\mathbb{E} \cong \Q He_{\Q}(\lambda)$ and $\Q He_{\Q}(\lambda)$ is isomorphic to $\Q(\zeta_{p})$ via the map \linebreak $xe_{\Q}(\lambda) \mapsto \zeta_{p},~ye_{\Q}(\lambda) \mapsto 1,~w^{2^{n-1}} \mapsto -1$. Also $\mathbb{F}=\mathcal{Z}(\Q Ge_{\Q}(\lambda^G))$ is isomorphic to the subfield $F=\Q (\zeta_{p} + \zeta_{p}^{r^t}+\cdots + \zeta_{p}^{r^{{(2^{n-1}-1)}t}})$ of  $\Q(\zeta_{p})$ kept fixed by $\langle \sigma \rangle$, where $\sigma:\Q(\zeta_{p}) \rightarrow \Q(\zeta_{p})$ is given by $\zeta_{p}$ to $\zeta_{p}^{r^t}$. This is because $\sigma_{w}$ sends $\operatorname{diag}(xe_{\Q}(\lambda),xe_{\Q}(\lambda),\cdots,xe_{\Q}(\lambda))_{p}$ to $z_{w}^{-1}\operatorname{diag}(xe_{\Q}(\lambda),xe_{\Q}(\lambda),\cdots,xe_{\Q}(\lambda))_{p}z_{w}$,\linebreak which is $\operatorname{diag}(x^{r^t}e_{\Q}(\lambda),x^{r^t}e_{\Q}(\lambda),\cdots,x^{r^t}e_{\Q}(\lambda))_{p}$ from the choice of $z_{w}=wA_{w}$. Hence, in view of equation (\ref{equ13}), it turns out that
$$(\mathbb{E}/\mathbb{F},\sigma_w, z_w^{2^{n-1}}) \cong  (\Q(\zeta_{p})/F, \langle \sigma\rangle,-p^{2^{n-2}}).$$
As $p \equiv 1 \operatorname{mod} 4,$  we have $\sqrt{p} \in \Q(\zeta_{p})$ and hence $p^{2^{n-2}} \in  N_{\Q(\zeta_{p})/F}(\Q(\zeta_{p})^{\times}),$ where $\Q(\zeta_{p})^{\times}$ are non zero elements of $\Q(\zeta_{p})$. Also $p$ is an odd prime and $[\mathbb{F}:\Q]$ is odd implies that $-1 \not\in  N_{\Q(\zeta_{p})/F}(\Q(\zeta_{p})^{\times})$. This gives $$ (\Q(\zeta_{p})/F, \langle \sigma\rangle,-p^{2^{n-2}}) \cong  (\Q(\zeta_{p})/F, \langle \sigma\rangle,-1)$$ and it is a central simple $F$-algebra  of Schur index 2, which, by (\cite{lam}, Page 226), is Brauer equivalent to the quaternion division algebra $\mathbb{H}(F)$ over $F$. Consequently,
$$\Q Ge_{\Q}(\lambda^G) \cong  M_{p2^{n-1}}(\mathbb{H}(F)).$$
This completes the proof for the Shoda pair $(H,K)=(Q,\langle y\rangle ).$ \para \noindent 
For the other shoda pair $(H,K)=(Q,\langle y,w^{2^{n-1}}\rangle ),$ if we consider  $\lambda : Q \to \mathbb{C}$ given by $x \mapsto \zeta_p, y \mapsto 1, w^{2^{n-1}} \mapsto 1$ and perform the similar computations, it turns out that $\Q Ge_{\Q}(\lambda^G) \cong  M_{p2^{n}}(F).$ \para \noindent
This section thus concludes with the following: 
\begin{theorem} If $G= G_{p^32^n}$, 
the Wedderburn decomposition of $\Q G$ is given by:
		$$\Q G \cong \oplus_{i=0}^{n} \Q(\zeta_{2^{n-i}}) \oplus (M_{2^n}(F) )^{\frac{p+1}{2}} \oplus M_{p2^{n-1}}(\mathbb{H}(F)) \oplus M_{p2^{n}}(F)$$
		where ${F}= \Q(\zeta_{p}+\zeta_{p}^{r^{t}}+\cdots+\zeta_{p}^{r^{(2^{n-1}-1)t}}),~ \mathbb{H}(F)$ is the quaternion algebra over $F$ and $(M_{2^n}(F) )^{\frac{p+1}{2}}$ is the direct sum of $\frac{p+1}{2}$ copies of $M_{2^n}(F).$
		\end{theorem}
\section{$G = (Q_{8} \times C_{7}) \rtimes C_{3}$}
Let $Q_{8}$ be the group of quaternions of order 8 and $C_{n}$ the cyclic group of order $n,~n \geq 1$. In this section, we consider $G$ to be the semidirect product of $Q_{8} \times C_{7}$ with $C_{3}$ given by the following presentation :
$$ G=\langle x,y,a,b~|~ x^4=y^4=b^7=a^3=e, x^2=y^2, xy^{-1}= yx, xb=bx, yb=by,$$ $$a^{-1}xa=y, a^{-1}ya=xy, a^{-1}ba=b^2\rangle .$$
Let $(H,K)=(\langle x,b\rangle ,1)$ and $\lambda$ a linear character on $H$ with kernel $K$ given by $x \mapsto \zeta_{4}, b \mapsto \zeta_{7}.$ Let $H_{1}=\langle x,y,b\rangle .$ Observe that $H \unlhd H_{1}$ and $I_{H_{1}}(\lambda)=H.$ Therefore, $\lambda^{H_{1}}$ is irreducible. Also, see that $H_{1} \unlhd G$ and $I_{G}(\lambda^{H_{1}}) = H_{1}$, which  implies that $(\lambda^{H_{1}})^{G} = \lambda^{G}$ is irreducible. Thus $(H,K)$ is a Shoda pair of $G.$ We see that $H \leq H_{1} \leq G$ is a strong inductive chain, in view of Remark \ref{remark2}$(i)$, as $H \unlhd H_{1}$ and $H_{1} \unlhd G.$ Let us compute $C_{0}=\operatorname{Cen}_{H_{1}}(e_{\Q}(\lambda))$ and $C_{1}=\operatorname{Cen}_{G}(e_{\Q}(\lambda^{H_{1}})).$ By (\cite{Jes}, Lemma 3.5.1 (3),(4)), $$C_{0}=\operatorname{Cen}_{H_{1}}(e_{\Q}(\lambda))= N_{H_{1}}(K) = H_{1}, ~\mbox{as}~ K=1.$$ It also implies that $e_{\Q}(\lambda^{H_{1}})= e_{\Q}(\lambda),$ i.e., $e_{\Q}(\lambda)$ is a primitive central idempotent of $\Q H_{1}.$  For any $g \in G,~e_{\Q}(\lambda)^g$ is thus a  primitive central idempotent of $\Q H_{1}^g.$ But $H_{1}^g = H_{1},$ as $H_{1} \unlhd G.$ So, $e_{\Q}(\lambda)^g$ is also a  primitive central idempotent of $\Q H_{1}$ for all $g \in G.$ Now if $g \in G \backslash \C_{G}(e_{\Q}(\lambda)),$ then $e_{\Q}(\lambda)$ and $e_{\Q}(\lambda)^g$ being distinct primitive central idempotents of $\Q H_{1}$ are orthogonal. This gives, in view of Lemma 3.5.2 of \cite{Jes}, that $\C_{G}(e_{\Q}(\lambda))= N_{G}(K) = G.$ Consequently, $$C_{1}= \C_{G}(e_{\Q}(\lambda^{H_{1}})) = \C_{G}(e_{\Q}(\lambda))=G.$$ Applying Theorems \ref{thm3} and \ref{thm5}, we hence obtain that
\begin{flalign*}
\Q Ge_{\Q}(\lambda^{G}) & \cong (\Q H_{0}e_{\Q}(\lambda) *^{\sigma_{H_{0}}}_{\tau_{H_{0}}} C_{0}/H_{0}) *^{\sigma_{H_{1}}}_{\tau_{H_{1}}} C_{1}/H_{1}, ~\mbox{as}~ k_{0}=k_{1}=1\\ & \cong \mathbb{E}*_{\tau} \G(\mathbb{E}/\mathbb{F}),
\end{flalign*} 
where $\mathbb{E}= \Q H_{0}e_{\Q}(\lambda),~\mathbb{F}=(\Q H_{0}e_{\Q}(\lambda)^{C_{0}/H_{0}})^{C_{1}/H_{1}}$ and $\G(\mathbb{E}/\mathbb{F})$ is of type $C_{0}/H_{0}$-by-$C_{1}/H_{1}$. As $C_{0}/H_{0}= H_{1}/H_{0} \cong C_{2}$ and $C_{1}/H_{1}=G/H_{1} \cong C_{3},$ it follows that $|\G(\mathbb{E}/\mathbb{F})| = 6.$ Also $\G(\mathbb{E}/\mathbb{F}),$ being a subquotient of the Galois group of a cyclotomic field extension of $\Q$, is abelian. Therefore, $\G(\mathbb{E}/\mathbb{F})$ is cyclic of order 6. \para \noindent We will now find a generator of $\G(\mathbb{E}/\mathbb{F}).$ For this purpose, we will find $A_{a} \in \Q H_{0}e_{\Q}(\lambda) * C_{0}/H_{0} (=\Q C_{0}e_{\Q}(\lambda))$ such that 
\begin{enumerate}[(i)]
	\item $aA_{a} \in \operatorname{Cen}_{\mathcal{A}}(\mathcal{B});$
	\item $z_{a}^{-1}\mathbb{E}z_{a}=\mathbb{E},$ where $z_{a}=aA_{a};$
	\item the automorphism $\sigma_{a} : \mathbb{E} \rightarrow \mathbb{E}$ given by $\alpha \mapsto z_{a}^{-1}\alpha z_{a}$ has order 6.
\end{enumerate} 
Once, we have found $A_{a}$ which satisfies (i)-(iii), then clearly $\sigma_{a}$ (as given in (iii) above) is a generator of $\G(\mathbb{E}/\mathbb{F}).$
\para \noindent In the present situation, as $k_{0}=k_{1}=1,$ we observe that $\mathcal{B}$ is just the center of $\mathcal{A}$ and hence $\operatorname{Cen}_{\mathcal{A}}(\mathcal{B})= \mathcal{A}.$ Therefore (i) holds for any arbitrary $A_{a} \in \Q H_{0}e_{\Q}(\lambda) * C_{0}/H_{0} (= \Q C_{0}e_{\Q}(\lambda))$.\para \noindent
In order to satisfy (ii) and (iii), it suffices to find $A_{a} \in \Q C_{0}e_{\Q}(\lambda)$ which satisfies the following:
\begin{equation}{\label{eq10}}
(aA_{a})^{-1}xe_{\Q}(\lambda)aA_{a} = x^3e_{\Q}(\lambda);\\
(aA_{a})^{-1}be_{\Q}(\lambda)aA_{a} = b^2e_{\Q}(\lambda).
\end{equation}
It can be seen that $A_{a}=xye_{\Q}(\lambda)+y^2e_{\Q}(\lambda)$ is one solution of the system of equations (\ref{eq10}). This gives us a generator of $\G(\mathbb{E}/\mathbb{F}),$ namely $\sigma_{a}$ which sends $\alpha$ in $\mathbb{E}$ to $z_{a}^{-1}\alpha z_{a}$ where $z_{a}=aA_{a}=a(xye_{\Q}(\lambda)+y^2e_{\Q}(\lambda)).$ \para \noindent Consequently, $\mathbb{E} *_{\tau} \G(\mathbb{E}/\mathbb{F})$ is isomorphic to the cyclic algebra $(\mathbb{E}/\mathbb{F},\sigma,z_{a}^6)$. A straight forward computation reveals that  $z_{a}^6= 8x^2e_{\Q}(\lambda).$ So,
\begin{equation}{\label{eq11}}
\Q Ge_{\Q}(\lambda^G) \cong (\mathbb{E}/\mathbb{F},\sigma_{a},8x^2e_{\Q}(\lambda)).
\end{equation}
Now $\mathbb{E}$ is isomorphic to $\Q(\zeta_{7}, \iota)$ via the map $xe_{\Q}(\lambda) \mapsto \iota,~be_{\Q}(\lambda) \mapsto \zeta_{7}$. Also, notice that $\sigma_{a}(xe_{\Q}(\lambda))= z_{a}^{-1}x e_{\Q}(\lambda)z_{a}=x^3e_{\Q}(\lambda)$ and $\sigma_{a}(be_{\Q}(\lambda))= z_{a}^{-1}be_{\Q}(\lambda)z_{a}=b^2e_{\Q}(\lambda).$ Therefore, 
\begin{equation}{\label{eq12}}
	(\mathbb{E}/\mathbb{F},\sigma_{a},8x^2e_{\Q}(\lambda)) \cong (\Q(\zeta_{7}, \iota)/\mathbb{F}), \sigma, -8)
\end{equation}
where $\sigma$ is the automorphism of $\Q(\zeta_{7}, \iota)$ which sends $\zeta_{7} \mapsto \zeta_{7}^2$ and $\iota \mapsto -\iota$ and $\mathbb{F}$ is the fixed field of $\Q(\zeta_{7}, \iota)$ by $\langle \sigma \rangle.$ \para \noindent
As $(\zeta_{7}+\zeta_{7}^2+\zeta_{7}^4)$ is fixed by $\sigma$, we have $\Q(\zeta_{7}+\zeta_{7}^2+\zeta_{7}^4) = \Q(\sqrt{-7})$ is contained in $\mathbb{F}.$ Also, $[\mathbb{F}:\mathbb{Q}]= \frac{[\mathbb{E}:\mathbb{Q}]}{[\mathbb{E}:\mathbb{F}]}= \frac{\phi([H:K])}{6}=\frac{\phi(28)}{6}=2$, which gives $\mathbb{F} \cong \Q(\sqrt{-7}).$ Here, $\phi$ is Euler's phi function. Since, $8 = \operatorname{N}_{\Q(\zeta_{7},\iota)/\Q(\sqrt{-7})}(\iota +1)$ and $-1 \notin \operatorname{N}_{\Q(\zeta_{7},\iota)/\Q(\sqrt{-7})}(\Q(\zeta_{7},\iota)^{\times}),$ it turns out that $(\Q(\zeta_{7},\iota)/\Q(\sqrt{-7}), \sigma, -8)$ is congruent to $(\Q(\zeta_{7},\iota)/\Q(\sqrt{-7}), \sigma, -1)$ has Schur index $2.$ Consequently, the dimension of the simple component $\Q Ge_{\Q}(\lambda^G)$ over its center being $|\G(\Q(\zeta_{7},\iota)/\Q(\sqrt{-7}))|^2,$ which is 36, we obtain that $$\Q Ge_{\Q}(\lambda^G) \cong (\Q(\zeta_{7},\iota)/\Q(\sqrt{-7}), \sigma, -1) \cong M_{3}(\mathbb{H}(\Q(\sqrt{-7}))).$$	 
\bibliographystyle{amsplain}
\bibliography{BG1}
\end{document}